\documentstyle{article}
\def\a{\alpha}            \def\b{\beta}
\def\g{\gamma}            \def\G{\Gamma}
\def\d{\delta}            \def\D{\Delta}
\def\th{\theta}           
\def\vth{\vartheta}       
\def\ep{\epsilon}         \def\vep{\varepsilon}
\def\vph{\varphi}
           
\def\m{\mu}

\def\lm{\lambda}          \def\L{\Lambda}

\def\tl{\tilde}           
\def\ul{\underline}

         \def\cH{{\cal H}}         
                  
         \def\cN{{\cal N}}         
                  \def\cR{{\cal R}}

\newcommand{\id}{\mathop{\rm id}\nolimits}
\newcommand{\ad}{\mathop{\rm ad}\nolimits}

\def\lb{\label}                 \def\ot{\otimes}
\def\cd{\cdots}                 \def\ld{\ldots}
\def\st{\stackrel}

             \def\ra{\rightarrow}
\def\lra{\leftrightarrow}

\def\ns{\normalsize}            \def\ls{\large}
                 
\def\trid{\bigtriangledown}
\def\ds{\displaystyle}
\def\uD{\underline{\Delta}}
\def\Dp{\Delta_{+}}
\def\uDp{\underline{\Delta}_{+}}
\def\vuDp{\vec{\underline{\Delta}}_{+}}

\def\Z2{Z\!\!\!\!\!\!Z_{\,2}}
\def\CC{{\bf C}\!\!\!{\rm l}\,}

\newcommand{\bn}{\begin{equation}}
\newcommand{\ed}{\end{equation}}
\newcommand{\bneqn}{\begin{eqnarray}}
\newcommand{\edeqn}{\end{eqnarray}}
\newcommand{\bnth}{\begin{theorem}}
\newcommand{\edth}{\end{theorem}}
\newcommand{\bnpr}{\begin{proposition}}
\newcommand{\edpr}{\end{proposition}}
\newcommand{\bnlm}{\begin{lemma}}
\newcommand{\edlm}{\end{lemma}}
\newcommand{\bndef}{\begin{definition}}
\newcommand{\eddef}{\end{definition}}
\newcommand{\bntb}{\begin{tabular}}
\newcommand{\edtb}{\end{tabular}}

\def\nin{\noindent}
\newtheorem{definition}{Definition}[section]
\newtheorem{lemma}{Lemma}[section]
\newtheorem{proposition}{Proposition}[section]
\newtheorem{theorem}{Theorem}[section]

\normalbaselines
\begin{document}

\nin
{\ls\bf PROJECTION OPERATOR METHOD
FOR\\[5pt]
QUANTUM GROUPS}\\[0.5cm]

\nin
\hangindent=2cm \hangafter=0
{V.N. TOLSTOY}\\
{\it Institute of Nuclear Physics\\
Moscow State University\\
119899 Moscow \& Russia\\}

\setcounter{equation}{0}
\section{Introduction}
At present there can be no doubt to say that the quantum group
theory has been one of the most important, modern and rapidly
developing directions of mathematics and mathematical physics in
the end of the twentieth century. Although initially (in the
early 1980's) the quantum group theory was formulated for solving
problems in the theory of integrable systems and statistical
physics, later the surprising connection of this theory with many
branches of mathematics, and the theoretical and mathematical
physics was discovered. Today the quantum group theory is
connected with such mathematical fields as special functions
(especially, with $q$-orthogonal polynomials and basic
hypergeometric series), the theory of difference and differential
equations, combinatorial analysis and representation theory,
matrix and operator algebras, noncommutative geometry, knot
theory, topology, category theory and so on. From point of view
of the mathematical physics there exists interconnection of the
quantum group theory with the quantum  inverse scattering method,
conformal and quantum fields theory and so on. It is expected
that the quantum groups will provide deeper understanding of
concept of symmetry in physics.

In the broad sense the notation "quantum groups" \cite{D} involves
different deformations of the universal enveloping algebras $U(g)$
of Lie algebras and superalgebras $g$, such as: $q$-deformations,
Yangians, elliptic and dynamic quantum groups, mixed deformations
(for example, two-parameter deformations) and so on.
In the restricted sense the quantum groups mainly mean the
$q$-deformations of the universal enveloping algebras $U_q(g)$ of Lie
algebras and superalgebras $g$ (sometimes the $q$-deformations of
the groups and supergroups are also included here).
In what follows we use the notation of the quantum groups in this
restricted sense.

It is well known that the method of projection operators for usual
(non-quantized) Lie algebras and superalgebras is powerful and
universal method for a solution of many problems in the
representation theory. For example, the method allows
to classify irreducible modules, to decompose modules on submodules
(e.g. to analyze structure of Verma modules), to describe reduced
(super)algebras (which are connected with reduction of a (super)algebra
to (super)subalgebra), to construct bases of modules (e.g. the
Gelfand-Tsetlin's type), to develop the detailed theory of Clebsch-Gordan
coefficients and another elements of Wigner-Racah calculus (including
compact analytic formulas of these elements and their symmetry properties)
and so on. It is evident that the projection operators of quantum
groups play the same role in their representation theory.

In these lectures we develop the projection operator method
for quantum groups. Here the term "quantum groups" means $q$-deformed
universal enveloping algebras of contragredient Lie (super)algebras of
finite growth (these (super)algebras include all finite-dimensional
simple Lie algebras and classical superalgebras, infinite-dimensional
affine (Kac-Moody) algebras and superalgebras).
Conventionally, contains of the lectures can be divided on two parts.
Basis fragments of the first part are: a combinatorial structure
of root systems, the $q$-analog of the Cartan-Weyl basis,
the extremal projector and the universal $R$-matrix for
any contragredient Lie (super)algebra of finite growth.
It should be noted that the explicit expressions for the extremal
projectors and the universal R-matrices are ordered products of special
$q$-series depending on noncommutative Cartan-Weyl generators.

In second part (Sects. 9-12) we consider some applications of the
extremal projectors. Here we use the projector operator method to
develop the theory of the Clebsch-Gordan coefficients for the
quantum algebras $U_q(su(2))$ and $U_q(su(3))$. In particular, we
give a very compact general formula for the canonical
$U_q(su(3))\supset U_q(su(2))$ Clebsch-Gordan coefficients in
terms of the $U_q(su(2))$ Wigner $3nj$-symbols which are connected
with the basic hyperheometric series. Then we apply the projection
operator method for the construction of the $q$-analog of the
Gelfand-Tsetlin basis for $U_q(su(n))$. Finally using analogy
between the extremal projector $p(U_q(sl(2)))$ of the quantum
algebra $U_q(sl(2))$ and the $\d(x)$-function we introduce
'adjoint extremal projectors' $p^{(n)}(U_q(sl(2))$
($n\!=\!1,2,\ld$) which are some generalizations of the extremal
projector $p(U_q(sl(2)))$, and which are analogies of the
derivatives of the $\d(x)$-function, $\d^{(n)}(x)$
($n\!=\!1,2,\ld$). The elements $p^{(n)}(U_q(sl(2))$ can be
applied to construction and description of decomposable
representations of quantum algebra $U_q(sl(2))$ (details see
\cite{DT}).

\setcounter{equation}{0}
\section{Preliminary information}
Let $g(A,\tau)$ be a contragredient Lie (super)algebra
of finite growth\footnote{These (super)algebras include all
finite-dimensional simple Lie algebras and classical superalgebras,
infinite-dimensional affine (Kac-Moody) algebras and superalgebras.}
with a symmetrizable Cartan matrix $A$ (i.e. $A\!=\!DA^{sym}$, where
$A^{sym}\!=\!(a_{\;ij}^{sym}\!)_{i,j\in I}$ is a symmetrical matrix,
and $D$ is an invertible diagonal matrix,
$D\!=\!{\rm diag}\,(d_1,d_2,\ld,d_r)$), $\tau\subset I$,
$I\!:=\!\{1,2,\ld,r\}$, and let $\Pi:=\{\a_{1},\ld,\a_{r}\}$
be a system of simple roots for $g(A,\tau)$.

The Lie (super)algebra $g:=g(A,\tau)$ and its universal enveloping
algebra $U(g)$ are completely determined by the Chevalley generators
$e_{\pm\a_{i}}'$, $h_{\a_{i}}'$ $(i =1,2,\ld,r)$ with the defining
relations \cite{K1}:
\bneqn
[h_{\a_{i}}',h_{\a_{j}}']&\!\!=\!\!&0~,\qquad\qquad\qquad\quad
[h_{\a_{i}}',e_{\pm\a_{j}}']=\pm a_{\;ij}^{sym}e_{\pm a_{i}}'~,
\lb{PC1}
\\[7pt]
[e_{\a_{i}}',e_{-\a_{j}}']&\!\!=\!\!&\d_{ij}h_{\a_{i}}',\qquad
(\ad e_{\pm\a_{i}}')^{n_{ij}}\,e_{\pm\a_{j}}'=0 \quad {\rm
for}\;\; i\neq j~, \lb{PC2} \edeqn where the positive integers
$n_{ij}$ are given as follows: $n_{ij}\!=\!1$ if
$a_{\,ii}^{sym}\!\!=\!a_{\,ij}^{sym}\!\!=\!0$, $n_{ij}\!=\!2$ if
$a_{\,ii}^{sym}\!\!=\!0,\;a_{\,ij}^{sym}\!\!\ne\!0$, and
$n_{ij}\!=\!-2a_{\,ij}^{sym}/a_{\,ii}^{sym}\!\!+\!1$ if
$a_{\,ii}^{sym}\!\!\neq\!0$. Moreover it is necessary also to add
all nontrivial relations of the form: \bn
[[e_{\pm\a_{i}}',e_{\pm\a_{j}}'],[e_{\pm\a_{i}}',e_{\pm\a_{k}}']]=0
\lb{PC3} \ed for each triple of simple roots $\a_{i}$, $\a_{j}$,
$\a_{k}^{}$ if they satisfy the condition: \bn
a_{\,ii}^{sym}\!=a_{\,jk}^{sym}\!=0~,\qquad
a_{\,ij}^{sym}\!=-a_{\,ik}^{sym}\!\ne0~. \lb{PC4} \ed Throughout
the paper the brackets $[\cdot,\cdot]$ and the symbol "$\ad$"
denote the supercommutator in $U(g)$, i.e. \bn (\ad
a)b\equiv[a,b]=ab-(-1)^{\deg(a)\deg(b)}ba \lb{PC5} \ed for all
homogeneous elements $a,\;b\in g$, where \bn
\deg(h_{\a_{i}}')\!=\!0\;\,{\rm for}\;i\in I,\quad
\deg(e_{\pm\a_i}')\!=\!0\;\,{\rm for}\;i\not\in\tau,\quad
\deg(e_{\pm\a_i}')\!=\!1\;\,{\rm for}\;i\in\tau. \lb{PC6} \ed
{\it Remarks}. (i) The triple relation (\ref{PC3}) may appear
only in the supercase for the following situation in the Dynkin
diagram:

\vskip -15pt
\bn
\setlength{\unitlength}{0.7mm}
\begin{picture}(100,15)(-50,0)
\put(-20,5){\makebox(0,0){$\alpha_{j}$}}
\put(0,5){\makebox(0,0){$\alpha_{i}$}}
\put(20,5){\makebox(0,0){$\alpha_{k}$}}
\put(-20,0){\circle*{0.8}}
\put(0,0){\makebox(0,0){$\otimes$}}
\put(20,0){\circle*{0.8}}
\put(-18,0){\makebox{\line (1,0){16}}}
\put(2,0){\makebox{\line (1,0){16}}}
\end{picture}
\label{PC7} \ed Here $\a_{i}$ is an odd gray root, and $\a_{j}$,
$a_{k}$ are not connected and they can be of any color (degree):
white, gray or dark, and moreover the lines connected the node
$\a_{i}$ with the nodes $\a_{j}$ and $a_{k}$ can be non-single.
The second equality (\ref{PC2}) is ordinary called the Serre
relation therefore the equality (\ref{PC3}) may be called the
triple Serre relation because it connects three root vectors
$e_{\a_i}$, $e_{\a_j}$ and $e_{\a_k}$.

\nin (ii) Besides the relations of type (\ref{PC3}) the
additional Serre relations of higher order can also occur but we
don't give them here because these relations appear only for the
Dynkin diagrams of special type. A total list of such diagrams
and corresponding additional Serre relations can be found in the
paper \cite{Y}.

Let $\Dp$ be a system of all positive roots of the (super)algebra
$g(A,\tau)$. Any root $\g$ of $\Dp$ has the form:
$\g=\sum_{i}^{r}n_i\a_i$, where all $n_i$ are nonnegative integers.
The total system of all roots, $\D$, has the form: $\D=\Dp\bigcup(-\Dp)$.
On the system $\D$ there is a bilinear form $(\cdot,\cdot)$ such that
$(\a_i,\a_j)\!=\!a_{ij}^{sym}$. The form is positive definite for
all simple finite-dimensional Lie algebras and it is nondegenerate
for all finite-dimensional contragredient Lie superalgebras.
With respect to this form the simple roots $\a_i\in\Pi$ are
classified (colored) as follows:
\begin{itemize}
\item
A simple root $\a_i$ is called even (white) if $(\a_i,\a_i)\ne0$ and
$2\a_i\notin\Dp$. (In this case $i\notin\tau$).
\item
A simple root $\a_i$ is called odd, dark if $(\a_i,\a_i)\ne0$ and
$2\a_i\in\Dp$. (In this case $i\in\tau$).
\item
A simple root $\a_i$ is called odd, grey if $(\a_i,\a_i)=0$.
One can show that doubled grey roots don't exist, $2\a_i\notin\Dp$.
(In this case $i\in\tau$).
\end{itemize}
The grey and dark roots occur only in the supercase.

Let $\g$ be any root of $\Dp$, this root is called odd if in its
decomposition on the simple roots, $\g=\sum_{i}^{r}n_{i}\a_{i}$,
the sum of the coefficients $n_i$ for all odd roots $\a_i$ is odd.
Otherwise the root $\g$ is called even. The parity of the negative
root $-\g$ coincides with the parity of the positive root $\g$.

Coloring of the roots is extended on all system $\D$ as follows:
\begin{itemize}
\item
All even roots are white. A white root is pictured by the white
node \mbox{\begin{picture}(10,10)\put(4,3){\circle{6}}\end{picture}}.
\item
An odd root $\g$ is called grey if $2\g$ is not any root.
This root is pictured by the grey node $\ot$.
\item
An odd root $\g$ is called dark if $2\g$ is a root.
This root is pictured by the dark node
\mbox{\begin{picture}(10,10)\put(4,3){\circle*{6}}\end{picture}}.
\end{itemize}
In the case of the affine Kac-Moody (super)algebras all roots
are also divided into real and imaginary. Every imaginary
root $n\d$ satisfies the condition $(n\d,\g)\!=\!0$ for all $\g\in\D$.
For the real roots this condition is not valid.

\setcounter{equation}{0}
\section{Combinatorial structure of root systems}
At first we remind the definition of the reduced system of the positive
root system $\D_+$ for any contragredient (super)algebras of finite
growth.
\bndef
The system $\ul{\D}_+$ is called the reduced system if it is defined
by the following way:
$\ul{\D}_+\!=\D_+\!\backslash\{2\g\in \D_+|\,\g\;{\rm is\; odd}\}$.
That is the reduced system $\ul{\D}_+$ is obtained from the total
system $\D_+$ by removing of all doubled roots $2\g$ where $\g$ is
a dark odd root.
\eddef
Combinatorial structure of root system of the contragredient
Lie (super)algebra of finite growth is connected with
notation of the normal ordering in the reduced system of positive roots.
\bndef
We say that the system $\underline{\Delta}_+$ is in normal
ordering if each composite (not simple) root
$\g\!=\!\alpha\!+\!\b$ ($\a,\b,\g\in\uDp$),
where $\alpha$ and $\b$ are not proportional roots ($\a\ne\lm\b$),
is written between its components $\alpha$ and $\b$.
It means that in the normal ordering system $\uDp$ we have either
\bn
\mbox{\small$\ld,\a,\ld,\a\!+\!\b,\ld,\b,\ld$}~,
\lb{CS1}
\ed
or
\bn
\mbox{\small$\ld,\b,\ld,\a\!+\!\b,\ld,\a,\ld$}\ .
\lb{CS2}
\ed
We say also that $\alpha\prec\b$ if $\alpha$ is located on the
left side of $\b$ in the normal ordering system $\uDp$,
i.e. this corresponds to the case (\ref{CS1}).
\eddef
The normal ordering system $\uDp$ is denoted by the symbol $\vuDp$.
It is evident that boundary (end) roots in $\vuDp$ are simple.
The combinatorial structure of the root system $\uDp$ is described
the following theorem.
\bnth
(i) Normal ordering in the system $\uDp$ exists for any mutual location
of the simple roots $\a_i,\; i=1,2,\ld,r$.

\nin
(ii) Any two normal orderings $\vuDp$ and $\vuDp'$ can be obtained one
from another by compositions of the following elementary inversions:
\bn
\mbox{\small$\a,\b\,\lra\,\b,\a$}~,
\lb{CS3}
\ed
\bn
\mbox{\small$\a,\a\!+\!\b,\b\,\lra\,\b,\a\!+\!\b,\a$}~,
\lb{CS4}
\ed
\bn
\mbox{\small$\a,\a\!+\!\b,\a\!+\!2\b,\b\,\lra\,
\b,\a\!+\!2\b,\a\!+\!\b,\a$}~,
\lb{CS5}
\ed
\bn
\mbox{\small$\a,\a\!+\!\b,2\a\!+\!3\b,\a\!+\!2\b,\a\!+\!3\b,\b\,\lra\,
\b,\a\!+\!3\b,\a\!+\!2\b,2\a\!+\!3\b,\a\!+\!\b,\a$},
\lb{CS6}
\ed
\bn
\begin{array}{l}
\mbox{\small$\a,\d\!+\!\a,2\d\!+\!\a,\ld,\infty\d\!+
\!\a,\d,2\d,3\d,\ld,\infty\d,
\infty\d\!-\!\a,\ld,2\d\!-\!\a,\d\!-\!\a\,\lra$}
\\[3pt]
\mbox{\small$\lra\,\d\!-\!\a,2\d\!-\!\a,\ld,\infty\d\!-
\!\a,\d,2\d,3\d,\ld,\infty\d,\infty\d\!+
\!\a,\ld,2\d\!+\!\a,\d\!+\!\a,\a$}~,
\end{array}
\lb{CS7}
\ed
\bn
\begin{array}{l}
\!\!\!
\mbox{\small$\a,\d\!+2\a,\d\!+\!\a,3\d\!+\!2\a,2\d\!+\!\a,\ld,
\infty\d\!+\!\a,(2\infty\!+\!1)\d\!+\!2\a,
(\infty\!+\!1)\d\!+\!\a,\d,2\d,\ld,$}
\\[5pt]
\!\!\!
\mbox{\small$\infty\d,(\infty\!+\!1)\d\!-\!\a,
(2\infty\!+\!1)\d\!-\!2\a,\infty\d-\a,\ld,
2\d\!-\!\a,3\d\!-\!2\a,\d\!-\!\a,\d\!-\!2\a,\,\lra$}
\\[5pt]
\!\!\!
\mbox{\small$\lra\d\!-\!2\a,\d\!-\!\a,3\d\!-\!2\a,2\d\!-\!\a,\ld,
\infty\d\!-\!\a,(2\infty\!+\!1)\d\!-\!2\a,(\infty\!+\!1)\d\!-\!\a,
\d,2\d,\ld,$}
\\[2pt]
\!\!\!
\mbox{\small$\infty\d,(\infty\!+\!1)\d\!+\!\a,
(2\infty\!+\!1)\d\!+\!2\a,\infty\d\!+\!\a,\ld,
2\d\!+\!\a,3\d\!+\!2\a,\d\!+\!\a,\d\!+\!2\a,\a$},
\end{array}
\lb{CS8}
\ed
where $\a\!-\!\b$ is not any root.
\edth
A proof of the second part {\it (ii)} for the case of the
finite-dimensional simple Lie algebras can be found in \cite{AST2}.
The full proof of the theorem is given in the outgoing paper
\cite{T5}.

The root systems in (\ref{CS3})--(\ref{CS8}) belong to the
(super)algebras of rank 2. The combinatorial theorem permits
to construct a $q$-analog of the Cartan-Weyl basis and to
reduce the proof of basic theorems for extremal projector and
the universal $R$-matrix for the quantum (super)algebra of
arbitrary rank to the proof of such theorems for the quantum
(super)algebras of rank 2.

\setcounter{equation}{0}
\section{Quantized Lie (super)algebras}
The quantum ($q$-deformed) universal enveloping (super)algebra
$U_{q}(g)$\footnote{We shall also use the name "the quantized Lie
(super)algebra $g$".}, where $g$ is a contragredient Lie
(super)algebra of finite growth, may be consider as a deformation
$f$ (reserving the grading) of the universal enveloping algebra
$U(g)$: $U(g)\st{f}{\mapsto}U_{q}(g)$
$(e_{\pm\a_{i}}'\st{f}{\mapsto}e_{\pm\a_{i}}^{}$,
$h_{\a_{i}}'\st{f}{\mapsto}h_{\a_{i}}^{})$, which modifies the
relations (\ref{PC2}) ¿ (\ref{PC3}). More precisely we have the
following definition \cite{T3,KT1}. \bndef The quantum
(super)algebra $U_{q}(g)$ (where $g:=g(A,\tau)$ is a
contragredient Lie (super)algebra of finite growth), is an
associative (super)algebra over $\CC[q,q^{-1}]$ with Chevalley
generators $e_{\pm\a_{i}}^{}$, $k_{\a_{i}}^{\pm 1}\!:=\!q^{\pm
h_{\a_{i}}}$, $(i\in I\!:=\!\{1,2,\ld,r\})$, and the defining
relations: 
\bneqn
k_{\a_{i}}^{}k_{\a_{i}}^{-1}&\!\!=\!\!&k_{\a_{i}}^{-1}k_{\a_{i}}^{}=1~,
\qquad\quad\; k_{\a_{i}} k_{\a_{j}}=k_{\a_{j}} k_{\a_{i}}~, \quad
\lb{QA1}
\\[7pt]
k_{\a_{i}}^{}e_{\pm\a_{j}}^{}k_{\a_{i}}^{-1}&\!\!=\!\!&
q^{\pm(\a_{i},\a_{j})}e_{\pm\a_{j}}^{}~,\quad
[e_{\a_{i}},e_{-\a_{j}}]=\d_{ij}
\mbox{\ls$\frac{k_{\a_{i}}^{}-k_{\a_{i}}^{-1}}{q-q^{-1}}$}~,
\lb{QA2}
\\[7pt]
(\ad_{q}e_{\pm\a_{i}})^{n_{ij}}e_{\pm\a_j}&\!\!=\!\!&0
\qquad\qquad\qquad\quad{\rm for}\;\;i\neq j~, \lb{QA3} \edeqn
where the positive integers $n_{ij}$ are the same as in the
relations (\ref{PC2}). Moreover, if any three simple roots
$\a_{i}$, $\a_{j}$, $\a_{k}^{}$ satisfy the condition (\ref{PC4})
then there are the additional triple relations of the form: \bn
[[e_{\pm\a_{i}},e_{\pm\a_{j}}]_q,[e_{\pm\a_{i}},e_{\pm\a_{k}}]_q]_q=0~.
\lb{QA4} \ed \eddef Here in (\ref{QA1})--(\ref{QA4}) the brackets
$[\cdot,\cdot]$ is the usual supercommutator (\ref{PC5}), and
$[\cdot,\cdot]_{q}$ and $\ad_{q}$ denote the $q$-deformed
supercommutator ($q$-supercommutator) in $U_{q}(g)$: \bn
(\ad_{q}e_{\a})e_{\b}\equiv[e_{\a},e_{\b}]_q=
e_{\a}e_{\b}-(-1)^{\deg(e_{\a})\deg(e_{\b})}q^{(\a,\b)}e_{\b}e_{\a}~,
\label{QA5} \ed where $(\a,\b)$ is a scalar product of the roots
$\a$ and $\b$, and the parity function $\deg(\cdot)$ is given by
\bn \deg(k_{\a_{i}})\!=\!0\;\,{\rm for}\;i\in I,\quad
\deg(e_{\pm\a_i})\!=\!0\;\,{\rm for}\;i\not\in\tau,\quad
\deg(e_{\pm\a_i})\!=\!1\;\,{\rm for}\;i\in\tau. \lb{QA6} \ed
Below we shall use the following short notation: \bn
\vth(\g):=\vth(e_{\g})=\deg(e_{\g})\ . \lb{QA7} \ed {\it
Remarks.} (i) It is not hard to verify that the the relations
(\ref{QA1})--(\ref{QA4}) are invariant with respect to the
replacement of $q$ by $q^{-1}$.

\nin
(ii) The outer $q$-supercommutator in (\ref{QA4}) is really the usual
supercommutator since $(\a_{i}\!+\!\a_{j}, \a_{i}\!+\!\a_{k})\!=\!0$.

\nin
(iii) The remark (ii) after the formula (\ref{PC7}) is also valid.

Clearly, the quantum (super)algebra $U_{q}(g)$ reduces to the usual
universal enveloping (super)algebra $U(g)$ if $q\ra 1$.

By direct calculations we can show that quantum (super)algebra
$U_{q}(g)$ is a Hopf (super)algebra with respect to a comultiplication
$\D_{q}$, an antipode $S_{q}$ and a counit $\vep$ defined as
\bn
\begin{array}{rcccl}
\D_{q}(k_{\a_i}^{\pm1})&\!\!=\!\!&k_{\a_i}^{\pm1}\ot k_{\a_i}^{\pm1}~,
\qquad\qquad\qquad\;
S_{q}(k_{\a_i}^{\pm1})&\!\!=\!\!&k_{\a_i}^{\mp1}~,
\\[5pt]
\D_{q}(e_{\a_i}^{})&\!\!=\!\!&
e_{\a_i}^{}\ot 1+k_{\a_i}^{-1}\ot e_{\a_i}^{}~,
\qquad\;\; S_{q}(e_{\a_i}^{})&\!\!=\!\!&-k_{\a_i}e_{\a_i}^{}~,
\\[5pt]
\D_{q}(e_{-\a_i}^{})&\!\!=\!\!&
e_{-\a_i}^{}\ot k_{\a_i}+1 \ot e_{-\a_i}^{}~,
\quad S_{q}(e_{-\a_i}^{})&\!\!=\!\!&-e_{-\a_i}^{}k_{\a_i}^{-1}~,
\lb{QA8}
\end{array}
\ed
\bn
\vep(e_{\pm\a_{i}})=0~,\qquad\vep(k_{\a_{i}})=\vep(1)=1~,
\phantom{possible}
\lb{QA9}
\ed

Both in the quantum and non-quantum case we can directly use the Chevalley
generators for construction of a 'monomial' basis in all universal
enveloping (super)algebra $U_{q}(g)$ ($U(g)$). Bases of such kind
were proposed by Verma for the non-quantized case and by Lusztig
\cite{L} for the general case.
The Lusztig basis is an universal one and it is called the
canonical basis. Both the Verma basis and the Lusztig basis have
rather complicated algebraic structure and therefore they
were not used in a broad fashion until now.
It is well known that a monomial basis constructed of Cartan-Weyl
generators is a more algebraically simple basis. Therefore a natural
problem is to construct a $q$-analog of the Cartan-Weyl basis
(the quantum Cartan-Weyl basis) for the quantum (super)algebra
$U_{q}(g)$.

Our method for construction of the $q$-analog of the Cartan-Weyl basis
and its general properties and also properties of the extremal projector
and the universal $R$-matrix are closely connected with the combinatorial
structure for the root system of the Lie (super)algebra $g$.

\setcounter{equation}{0}
\section{Quantum Cartan-Weyl basis}
The $q$-analog of the Cartan-Weyl basis for $U_{q}(g)$ is constructed
by using the following inductive algorithm \cite{T3},
\cite{KT1}--\cite{KT5}.

\nin
{\it We fix some normal ordering $\vuDp$ and put by induction
\bn
e_{\g}:=[e_{\a},e_{\b}]_{q},\qquad\quad
e_{-\g}:=[e_{-\b},e_{-\a}]_{{q}^{-1}}
\lb{CW1}
\ed
if $\g\!=\!\a\!+\!\b$, $\a\prec\g\prec\b$ ($\a,\b,\g\in\uDp$), and
the segment $[\alpha;\b]\subseteq\vuDp$ is minimal one including
the root $\g$, i.e. the segment has not another roots $\alpha'$ and $\b'$
such that $\a'\!+\!\b'\!=\!\g$. Moreover we put
\bn
k_{\g}:=\prod_{i=1}^{r} k_{\a_{i}}^{l_{i}}~,
\label{CW2}
\ed
if $\g=\sum_{i=1}^{r}l_{i}\a_{i}$ ($\g\in\uDp,\;\a_i\in \Pi$)}.

\noindent By this procedure one can construct the total quantum
Cartan-Weyl basis for all quantized finite-dimensional simple
contragredient Lie (super)algebras. In the case of the quantized
infinite-dimensional affine Kac-Moody (super)algebras we have to
apply one more additional condition. Namely, first we construct
all root vectors $e_\g$ ($\g\in\uD$) by means of the given
procedure, and then we overdeterminate the generators $e_{n\d}$
of the imaginary roots $n\d\in\uD$ in a way that the new
generators $e_{n\d}'$ are mutually commutative if they are not
conjugate generators. Because of the fact that we do not have a
sufficient place here to describe the overdetermination of
imaginary root generators in details, we are restricted to a
consideration of finite-dimensional case, i.e. when $g$ is a
finite-dimensional simple contragredient Lie (super)algebra.

The quantum Cartan-Weyl basis is characterized by the following
properties \cite{KT1}--\cite{KT5}.
\begin{proposition}
The root vectors $\{e_{\pm\g}\}$ ($\g\,\in\uDp$) satisfy the following
relations:
\bneqn
k_{\a}^{\pm 1}e_{\g}&=&q^{\pm(\a,\g)}e_{\g}
k_{\a}^{\pm 1}~,
\lb{CW3}
\\[7pt]
[e_{\g}^{},e_{-\g}^{}]&=&a(\g)\frac{k_{\g}-k_{\g}^{-1}}
{q-q^{-1}}~,
\lb{CW4}
\edeqn
\bn
{[e_{\a}^{},e_{\b}^{}]}_{q}=\sum_{\a\prec\g_{1}\prec\ld
\prec\g_{n}\prec\b}C_{m_{i},\g_{i}}e_{\g_{1}}^{m_{1}}
e_{\g_{2}}^{m_{2}}\cd e_{\g_{n}}^{m_{n}},
\lb{CW5}
\ed
where $\sum_{i}^{n}m_{i}\g_{i}\!=\!\a\!+\!\b$, and
the coefficients $C_{\cd}$ are rational functions of $q$ and they
do not depend on the Cartan elements $k_{\a_{i}},\: i=1,2,\ld n$,
and also
\bn
{[e_{\b },e_{-\a}]}=
\sum C_{m_{i},\g_{i};m_{j}',\g_{j}'}'e_{-\g_{1}}^{m_{1}}
e_{-\g_{2}}^{m_{2}}\cd e_{-\g_{p}}^{m_{p}}
e_{\g_{1}'}^{m_{1}'}e_{\g_{2}'}^{m_{2}'}\cd
e_{\g_{s}'}^{m_{s}'}
\lb{CW6}
\ed
where the sum is taken on $\g_{1},\ld,\g_{p},\g_{1}',\ld,\g_{s}'$ and
$m_{1},\ld,m_{p}, m_{1}',\ld,m_{s}'$ such that
\[
\g_{1}\prec\ld\prec\g_{p}\prec\a\prec\b\prec\g_{1}'\prec\ld\prec\g_{s}'~,
\quad\sum_{l}(m_{l}'\g_{l}'-m_{l}\g_{l})=\b-\a
\]
and the coefficients $C_{\cd}'$ are rational functions of $q$ and
$k_{\a}$ or $k_{\b}$. The monomials $e_{\g_{1}}^{n_{1}}
e_{\g_{2}}^{n_{2}} \cd e_{\g_{p}}^{n_{p}}$ and
$e_{-\g_{1}}^{n_{1}}e_{-\g_{2}}^{n_{2}}\cd
e_{-\g_{p}}^{n_{p}}$, ($\g_{1} \prec\g_{2}\prec\cd
\prec\g_{p}$), generate (as a linear space over $U_q(\cH)$)
subalgebras $U_{q}(b_{+})$ and $U_{q}(b_{-})$ correspondingly.
The monomials
\[
e_{-\g_{1}}^{n_{1}}e_{-\g_{2}}^{n_{2}}\cd
e_{-\g_{p}}^{n_{p}}e_{{\g'}_{1}}^{{n'}_{1}}
e_{{\g'}_{2}}^{{n'}_{2}}\cd e_{{\g'}_{s}}^{{n'}_{s}},
\]
where $\g_{1}\prec\g_{2}\prec\cd\prec\g_{p}$ and
$\g'_{1} \prec\g'_{2}\prec\cd\prec\g'_{s}$), generate
$U_{q}(g)$ over $U_{q}(\cH)$.
\end{proposition}
Here the algebra $U_q(\cH)$ is generated by the Cartan elements
$k_{\a_i}$($i\!=\!1,2,\ld,r$).

Now we consider some extensions of $U_{q}(g)$,
$U_{q}(b_{+})\ot U_{q}(b_{-})$ and $U_{q}(g)\ot U_{q}(g)$ since
the extremal projector and the universal $R$-matrix are elements of these
extensions.

\setcounter{equation}{0}
\section{Taylor extensions of $U_{q}(g)$, $U_{q}(b_{+})\ot
U_{q}(b_{-})$ and $U_{q}(g)\ot U_{q}(g)$}

Let $\mbox{Fract}\,(U_{q}(K))$ be a field of fractions
over $U_{q}(K)$, i.e. $\mbox{Fract}\,(U_{q}(K))$ is
an associative algebra of rational functions of the elements
$k_{\a_{i}}^{\pm 1}$, $(i=1,2,\ld,r)$. We put
\bn
\tl{U}_{q}(g)=\mbox{Fract}\,(U_{q}(K))\ot_{U_{q}(K)}U_{q}(g)~.
\lb{TE1}
\ed
Evidently, the extension $\tl{U}_{q}(g)$ is an associative algebra.
The algebra $\tl{U}_{q}(g)$ is called the Cartan extension of
the quantum algebra $U_{q}(g)$.

Let $\{e_{\pm\g}\}$, $\g\in\uD_{+}$, be the root vectors of
the quantum Cartan-Weyl basis built in accordance with some fixed
normal ordering in $\uD_{+}$. Let us construct a formal Taylor series
on the following monomials
\bn
e_{-\b}^{n_{\b}}\cd e_{-\g}^{n_{\g}}
e_{-\a}^{n_{\a}}\:e_{\a}^{m_{\a}}
e_{\g}^{m_{\g}}\cd e_{\b}^{m_{\b}}
\lb{TE2}
\ed
with coefficients from $\mbox{Fract}\,(U_{q}(K))$, where
$\a\prec\g\prec\cd\prec\b$ in a sense of the fixed normal ordering in
$\uD_{+}$ and nonnegative integers
$n_{\b},n_\g\ld,n_{\a},m_{\a},m_\g\ld,m_{\b}$
are subjected to the constraints
\bn
\Bigr|\sum_{\g\in\uD_{+}}(n_{\g}-m_{\g})c_{i}^{(\g)}
\Bigl|\leq\mbox{const}~,\qquad i=1,2,\cd,r~,
\label{TE3}
\ed
where $c_{i}^{(\g)}$ are coefficients in a decomposition of the root
$\g$ with respect to the system  of simple roots $\Pi$.
Let $T_{q}(g)$ be a linear space of all such formal series.
We have the following simple proposition.
\begin{proposition}
The linear space $T_{q}(g)$ is an associative algebra with respect to
a multiplication of formal series.
\end{proposition}
The algebra $T_{q}(g)$ is called the Taylor extension of $U_{q}(g)$.

Let $\mbox{Fract}\, (U_{q}(K\ot K))$ be a field of fractions
generated by the following elements: $1\ot k_{\a_{i}}$,
$k_{\a_{i}}\ot 1$ and $q^{h_{\a_{i}} \ot h_{\a_{j}}}$,
($i,j=1,2,\ld,r$). Let us consider a formal Taylor series of the
following monomials \bn e_{\a}^{n_{\a}}e_{\g}^{n_{\g}}\cd
e_{\b}^{n_{\b}}\ot e_{-\b}^{m_{\b}}\cd
e_{-\g}^{m_{\g}}e_{-\a}^{m_{\a}} \lb{TE4} \ed with coefficients
from $\mbox{Fract}\,(U_{q}(K\ot K))$, where
$\a\prec\g\prec\cd\prec\b$ in a sense of the fixed normal ordering
in $\uDp$ and nonnegative integers $n_{\b},\ld,n_{\a},m_{\a},
\ld,m_{\b}$ are subjected to the constraint (\ref{TE3}). Let
$T_{q}(b_{+}\ot b_{-})$ be a linear space of all such formal
series. The following proposition holds.
\begin{proposition}
The linear space $T_{q}(b_{+}\ot b_{-})$ is an associative algebra
with respect to a multiplication of formal series.
\end{proposition}
The algebra $T_{q}(b_{+}\ot b_{-})$ will be called the Taylor extension
of $U_{q}(b_{+})\ot U_{q}(b_{-})$.

At least we consider a formal Taylor series of the following monomials
\bn
e_{-\b}^{m_{\b}}\cd e_{-\g}^{m_{\g}}e_{-\a}^{m_{\a}}e_{\a}^{n_{\a}}
e_{\g}^{n_{\g}}\cd e_{\b}^{n_{\b}}\ot e_{-\b}^{m_{\b}'}\cd e_{-\g}^{m_{\g}'}
e_{-\a}^{m_{\a}'}e_{\a}^{n_{\a}'}e_{\g}^{n_{\g}'}\cd e_{\b}^{n_{\b}'}
\lb{TE5}
\ed
with coefficients from $\mbox{Fract}\,(U_{q}(K\ot K))$, where
$\a\prec\g\prec\cd\prec\b$ in a sense of the fixed normal ordering in
$\Delta_{+}$ and nonnegative integers $n_{\b},\ld,n_{\a},
m_{\a},\ld,m_{\b}$ and $n_{\b}',\ld,n_{\a}', m_{\a}',\ld,m_{\b}'$
are subjected to the constraints
\bn
\Bigl|\sum_{\g\in\Delta_{+}}(n_{\g}+n_{\g}'-m_{\g}-
m_{\g}')c_{i}^{(\g)}\Bigr|\leq\mbox{const}~,\qquad i=1,2,\cd,r~.
\lb{TE6}
\ed
Let $T_{q}(g\ot g)$ be a linear space of all such formal series.
The following simple proposition holds.
\begin{proposition}
The linear space $T_{q}(g \ot g)$ is an associative algebra
with respect to a multiplication of formal series.
\end{proposition}
The algebra $T_{q}(g\ot g)$ will be called the Taylor extension
of $U_{q}(g)\ot U_{q}(g)$.
Evidently the following embedding hold
\bn
\begin{array}{ccc}
&T_{q}(g\ot g)\supset T_{q}(b_{+}\ot b_{-})~,&
\\[7pt]
&T_{q}(g\ot g)\supset T_{q}(g)\ot T_{q}(g)\supset\D_{q}(T_{q}(g))~.&
\end{array}
\lb{TE7}
\ed

\setcounter{equation}{0}
\section{Extremal projector}
By definition, the extremal projector for $U_{q}(g)$ is a nonzero
element $p\!:=\!p(U_q(g))$ of the Taylor extension $T_q(g)$,
satisfying the equations \bn
e_{\a_{i}}p=pe_{-\a_{i}}=0\quad(\forall\;\a_{i}\in\Pi)~, \qquad
p^{2}=p~. \lb{EP1} \ed Acting by the extremal projector $p$ on
any highest weight $U_{q}(g)$-module $M$ we obtain a space
$M^{0}=pM$ of highest weight vectors for $M$ (if $pM$ has no
singularities).

Fix some normal ordering $\vuDp$ and let $\{e_{\pm\g}\}$ ($\g\in\uDp$)
be the corresponding Cartan-Weyl generators. The following statement
holds for any quantized finite-dimensional contragredient Lie
(super)algebra\footnote{The theorem is also valid for the quantized
infinite-dimensional affine Kac-Moody (super)algebras, but in this case
the formulas (\ref{EP3}) and (\ref{EP4}) for the imaginary roots
$\g=n\d$ should be more detailed (see \cite{KLT,T5} as examples).}
$g$ \cite{T3,KT2}.

\begin{theorem}
The equations (\ref{EP1}) have a unique nonzero solution in the space
of the Taylor extension $T_{q}(g)$ and this solution  has the form
\bn
p=\!\!\prod_{\g\in\vuDp}\!\!p_{\g}^{}~,
\lb{EP2}
\ed
where the order in the product coincides with the chosen normal
ordering of $\uDp$ and the elements $p_{\g}$ are
defined by the formulae
\bn
p_{\g}^{}=\sum_{m\geq0}\mbox{\ls$\frac{(-1)^{m}}
{(m)_{\bar{q}_{\g}^{}}!}$}\,\vph_{\g,m}^{}e_{-\g}^{\,m}e_{\g}^{m}~,
\lb{EP3}
\ed
\bneqn
\vph_{\g,m}^{}&\!\!=\!\!&\mbox{\ls$\frac{(q-q^{-1})^{m}
q^{-\frac{1}{4}m(m-3)(\g,\g)}q^{-m(\rho,\g)}}
{(a(\g))^m\prod\limits_{l=1}^{m}\Bigl(k_{\g}^{}
q^{(\rho,\g)+\frac{l}{2}(\g,\g)}
-(-1)^{(l-1)\th(\g)}k_{\g}^{-1}
q^{-(\rho,\g)-\frac{l}{2}(\g,\g)}\Bigr)}$}~.
\lb{EP4}
\edeqn
Here $\rho$ is a linear function such that
$(\rho,\a_{i})\!=\!\frac{1}{2}(\a_{i},\a_{i})$ for all
simple roots $\a_i\in\Pi$; $a(\g)$ is a factor in the relation
(\ref{CW4}); $\bar{q}_{\g}^{}\!:=\!(-1)^{\th(\g)}q^{-(\g,\g)}$;
the symbol $(m)_q$ is given by the formula:
\bn
(m)_{q}:=\mbox{\ls$\frac{q^{m}-1}{q-1}$}~.
\lb{EP5}
\ed
\end{theorem}
In the limit $q\to 1$ we obtain the extremal projector for the
(super)algebra $g$: $\lim\limits_{q\to1}p(U_q(g))=p(g)$
\cite{AST1,AST2,T1,T2}.
%
A proof of the theorem actually reduces to the proof for the case
of the quantized (super)algebras of rank 2, and it is similar to the
case of non-deformed finite-dimen\-sional simple Lie algebras
\cite{AST2}.

\setcounter{equation}{0}
\section{Universal $R$-matrix}
By definition, the universal $R$-matrix for the Hopf (super)algebra
$U_{q}(g)$ is an invertible element of the Taylor extension
$T_{q}(b_{+}\ot b_{-})$, satisfying the equations
\bn
\tl{\D}_q(x)=R\D(x)R^{-1},\qquad\forall\;x\in U_{q}(g),
\lb{UR1}
\ed
\bn
(\D_q\ot\id)R=R^{13}R^{23}~,\qquad (\id\ot\D_q)R=R^{13}R^{12},
\lb{UR2}
\ed
where $\tl{\D}_q$ is an opposite comultiplication:
$\tl{\D}_q=\sigma\D_q$, $\sigma(x\ot y)=
(-1)^{\deg x \deg y} y\ot x$ for all homogeneous elements
$x,y \in U_{q}(g)$. In (\ref{UR2}) we use standard notation
$R^{12}\!=\!\sum a_{i}\ot b_{i}\ot1$,
$R^{13}\!=\!\sum a_{i}\ot1\ot b_{i}$,
$R^{23}\!=\!\sum 1\ot a_{i}\ot b_{i}$
if $R$ has a form $R=\sum a_{i}\ot b_{i}$.

\nin
We employ the following standard notation for the q-exponential:
\bn
\exp_{q}(x):=1+x+\mbox{\ls$\frac{x^{2}}{(2)_{q}!}$}+\ld+
\mbox{\ls$\frac{x^{n}}{(n)_{q}!}$}+\ld=\sum_{n\geq0}
\mbox{\ls$\frac{x^{n}}{(n)_{q}!}$}~,
\lb{UR3}
\ed
where $(n)_{q}$ is defined by the formulas (\ref{EP5}).

Fix some normal ordering $\vuDp$ and let $\{e_{\pm\g}\}$ ($\g\in\uDp$)
be the corresponding Cartan-Weyl generators. The following statement
holds for any quantized finite-dimensional contragredient Lie
(super)algebra\footnote{The theorem is also valid for the quantized
infinite-dimensional affine Kac-Moody (super)algebras, but in
this case the formula (\ref{UR5}) for the imaginary roots $\g=n\d$
should be more detailed (see \cite{TK,KT4,KT5}.} $g$ \cite{KT1,KT2,KT3}.
\begin{theorem}
The equation (\ref{UR1}) has a unique (up to a multiplicative
constant) invertible solution in the space of the Taylor extension
$T_{q}(b_{+}\ot b_{-})$ and this solution has the form \bn
R=\Bigr(\prod_{\g\in\vuDp}R_{\g}\Bigr)\cdot K, \lb{UR4} \ed where
the order in the product coincides with the chosen normal
ordering $\vuDp$ and the elements $R_{\g}$ and $K$ are defined by
the formulas 
\bn 
R_{\g}=
\exp_{\bar{q}_{\g}^{}}\left((-1)^{\th(\g)}(q-q^{-1})
(a(\g))^{-1}(e_{\g}\ot e_{-\g})\right),
\lb{UR5} 
\ed 
\bn
K=q^{\sum_{i,j}d_{ij}(h_{\a_{i}}\ot h_{\a_{j}})} 
\lb{UR6} 
\ed
where $a(\g)$ is a factor from the relation (\ref{CW4}), and
$d_{ij}$ is an inverse matrix for a symmetrical Cartan matrix
$(a_{ij}^{sym})$ if $(a_{ij}^{sym})$ is not degenerated. (In a
case of a degenerated $(a_{ij}^{sym})$ we extend it up to a
non-degenerated matrix $(\tilde{a}_{ij}^{sym})$ and take an
inverse to this extended  matrix). Moreover the solution
(\ref{UR4}) is the universal $R$-matrix, i.e. it satisfies the
equations 
(\ref{UR2}) too.
\end{theorem}
A proof of the theorem actually reduces to the proof for the case
of the (super)algebras of rank 2 (see \cite{KT1}).

In the rest sections we consider some applications of the extremal
projectors.

\setcounter{equation}{0}
\section{Clebsch-Gordan and Racah coefficients for
the quantum algebras $U_q(su(2))$}
Let $J_{\pm}$, $q^{\pm J_{0}}$ be generators of the quantum algebra
$U_q(su(2))$. These generators satisfy the standard relations:
\bn
\begin{array}{rcl}
q^{J_0}J_{\pm}&\!\!=&\!\!q^{\pm1}J_{\pm}q^{J_0},\quad\;\;
[J_{+},J_{-}]=\frac{q^{2J_0}-q^{-2J_0}}{q-q^{-1}}\equiv [2J_0]~,
\\[7pt]
J_{\pm}^*&\!\!=&\!\!J_{\mp}~,\qquad\quad\;\;\; J_0^*=J_0~,
\qquad\quad\;\; q^*=q\;({\rm or}\;\,q^{-1})~.
\lb{su21}
\end{array}
\ed
Here and in what follows we use the notation:
$[x]\!=\!(q^x\!-\!q^{-x})/(q\!-\!q^{-1})$.
The Hopf structure of $U_q(su(2))$ is given by the following formulas
for the comultiplication $\D_q$, and the antipode $S_q$:
\bn
\begin{array}{rcccl}
\D_{q}(J_0)&\!\!=\!\!&J_0\ot 1+ 1\ot J_0~,
\qquad\qquad\quad\;
S_{q}(J_0)&\!\!=\!\!&-J_0~,
\\[5pt]
\D_{q}(J_{\pm}^{})&\!\!=\!\!&
J_{\pm}^{}\ot q^{J_0}+q^{-J_0}\ot J_{\pm}^{}~,\qquad\;
S_{q}(J_{\pm}^{})&\!\!=\!\!&-q^{\pm 1}J_{\pm}^{}~.
\lb{su22}
\end{array}
\ed
Let $\{|jm\bigr>\}$ be a canonical basis of the $U_q(su(2))$-irreducible
representation (IR) with the spin $j$. These basis vectors satisfy the
relations:
\bn
\begin{array}{rcl}
q^{J_0}|jm\bigr>&\!\!=\!\!&q^{m}|jm\bigr>~,
\\[7pt]
J_{\pm}|jm\bigr>&\!\!=\!\!&\sqrt{[j\mp m][j\pm m+1]}|jm\pm1\bigr>~.
\lb{su23}
\end{array}
\ed The vector $|jm\bigr>$ can be represented in the following
form 
\bn 
|jm\bigr>=F_{\!m;j}^{\,j}|jj\bigr>~, 
\lb{su24} 
\ed 
where
\bn
F_{\!m;j}^{\,j}=
\mbox{\ls$\sqrt{\frac{[j+m]!}{[2j]![j-m]!}}$}\,J_{-}^{j-m}~,
\lb{su25} 
\ed
and $|jj\bigr>$ is the highest weight vector, i.e.
\bn 
J_{+}^{}|jj\bigr>=0~. 
\lb{su26} 
\ed 
The operator
$F_{\!m;j}^{\,j}$ is called the lowering operator. We can also
introduce the rising operator 
\bn
F_{\!j;m}^{\,j}=
\mbox{\ls$\sqrt{\frac{[j+m]!}{[2j]![j-m]!}}$}\,J_{+}^{j-m}~,
\lb{su27}
\ed 
which has the property 
\bn
|jj\bigr>=F_{\!j;m}^{\,j}|jm\bigr>~. 
\lb{su28} 
\ed 
The extremal
projector $p$ for $U_q(su(2))$ can be represented in the form 
\bn
p=\sum_{n=0}^{\infty}\mbox{\ls$\frac{(-1)^n\bar{\G}_q(2J_0+2)}
{[n]!\bar{\G}_q(2J_0+n+2)}$}\,J_{-}^nJ_{+}^n~, 
\lb{su29} 
\ed 
where
$\bar{\G}_q(x)$ is the modified $q$-gamma function \bn
\bar{\G}_q(x+1)=[x]\G_q(x)~. \lb{su210} \ed This function is
connected with the standard Heine-Thomae $q$-gamma function
$\G_q(x)$ by the relation $\bar{\G}_q(x)=q^{x(x-1)/{4}}\G_q(x)$.
The extremal projector $p$ satisfies the relations: 
\bn
J_{+}^{}p=pJ_{-}^{}=0~,\qquad\quad p^{2}=p~. 
\lb{su211} 
\ed 
We multiply the extremal projector $p$ by the lowering an rising
operators as follows \bn
P_{\!m;m'}^{\,j}:=F_{\!m;j}^{\,j}\,p\,F_{\!j;m'}^{\,j}~. 
\lb{su212} 
\ed
Below we assume that the operator $P_{\!m;m'}^{\,j}$ acts in a
vector space of the weight $m'$. The operator $P_{\!m;m'}^{\,j}$
is called the general projection operator.

Let $\{|j_im_i\bigr>\}$ be canonical bases of two IRs $j_i$
$(i=1,2)$. Then $\{|j_1m_1\bigr>|j_2m_2\bigr>\}$ be an 'uncoupled'
bases in the representation $j_1\ot j_2$ of $U_q(su(2))\ot
U_q(su(2))$. In this representation there is another basis
$|j_1j_2\!:\!j_3m_3\bigr>_{\!q}$ which is called a 'coupled' basis 
with respect to $\D_q(U_q(su(2)))$. We can expand the coupled basis 
in terms of the uncoupled basis $\{|j_1m_1\bigr>|j_2m_2\bigr>\}$: 
\bn
\bigl|j_1j_2\!:\!j_3m_3\bigr>_{\!q}=\sum_{m_1,m_2}^{}
\bigl(j_1m_1\,j_2m_2|j_3m_3\bigr)_{\!q}\bigr|j_1m_1\bigr>
\bigr|j_2m_2\bigr>~, 
\lb{su213} 
\ed 
where the matrix element
$\bigl(j_1m_1\,j_2m_2|j_3m_3\bigr)_{\!q}$ is called the
Clebsch-Gordan coefficient (CGC). After some manipulations we can
show that CGC is presented by 
\bn
\bigl(j_1m_1\,j_2m_2\bigr|j_3m_3\bigr)_{\!q}=
\mbox{\ls$\frac{\langle j_1m_1|\langle
j_2m_2|\D_q(P_{\!m_3;j_3}^{j_3})|j_1j_1\rangle|j_2j_3-j_1\rangle}
{\sqrt{\langle j_1j_1|\langle
j_2\,j_3-j_1|\D_q(P_{\!j_3;j_3}^{j_3})
|j_1j_1\rangle|j_2\,j_3-j_1\rangle}}$}~. 
\lb{su214} 
\ed 
This is a formula for calculation of CGCs. Using the explicit 
expression (\ref{su212}) for the general projection operator
$P_{m_3;j_3}^{j_3}$, the formulas (\ref{su22}) for the
comultiplication $\D_q$ and the actions (\ref{su23}) for the
generators of $U_q(su(2))$ on the canonical basis vectors
$\big|j_im_i\bigr>$ $(i\!=\!1,2)$ it is not hard to calculate the
numerator and the denominator of the right side of (\ref{su214}).
As result we obtain the following expression for CGC of the
quantum algebra $U_q(su(2))$: 
\bn
\begin{array}{rcl}
&&\!\!\!\!\!\!\!
\bigl(j_1m_1\,j_2m_2\bigr|j_3m_3\bigr)_{\!q}=\delta_{m_1\!+m_2,m_3}^{}
q^{-\frac{1}{2}(j_1\!+\!j_2\!-\!j_3)
(j_1\!+\!j_2\!+\!j_3\!+\!1)+\!j_1m_2\!-\!j_2m_1}
\\[7pt]
&&
\phantom{possi}
\times\mbox{\ls$\sqrt{\frac{[2j_3+1][j_1+j_2-j_3]![j_1-j_2+j_3]!
[j_1+j_2+j_3+1]![j_2-m_2]![j_3+m_3]!}
{[-j_1+j_2+j_3]![j_1+m_1]![j_1-m_1]![j_2+m_2]![j_3-m_3]!}}$}
\\[7pt]
&&
\phantom{possi}
\times\sum\limits_{n}\mbox{\ls$\frac{(-1)^{j_1\!+\!j_2\!-\!j_3\!-\!n}
q^{n(j_1\!+m_1)}[2j_2-n]![j_1+j_2-m_3-n]!}
{[n]![j_1+j_2-j_3-n]![j_2-m_2-n]![j_1+j_2+j_3+1-n]!}$}~.
\lb{su215}
\end{array}
\ed
A total list of different explicit expressions and symmetry properties
for the $q$-CGCs can be found, for example, in \cite{KR,STK1,STK2,STK3}.

The general formula (\ref{su215}) can be expressed in terms of
the basic hypergeometric series
\bn
\begin{array}{rcl}
&&\!\!\!\!\!\!\!
\bigl(j_1m_1\,j_2m_2\bigr|j_3m_3\bigr)_{\!q}=\delta_{m_1\!+m_2,m_3}^{}
q^{-\frac{1}{2}(j_1\!+\!j_2\!-\!j_3)
(j_1\!+\!j_2\!+\!j_3\!+\!1)+\!j_1m_2\!-\!j_2m_1}
\\[7pt]
&&
\phantom{possibleaa}
\times(-1)^{j_1\!+j_2\!-j_3}\mbox{\ls$
\sqrt{\frac{[2j_3+1][j_1-j_2+j_3]!([2j_2]!)^2}
{[j_1+j_2-j_3]![-j_1+j_2+j_3]![j_1+j_2+j_3+1]!}}$}
\\[7pt]
&&
\phantom{possibleaa}
\times\mbox{\ls$\sqrt{\frac{[j_3+m_3]!([j_1+j_2-m_3]!)^2}
{[j_1+m_1]![j_1-m_1]![j_2+m_2]![j_2-m_2]![j_3-m_3]!}}$}
\\[7pt]
&&
\phantom{possibleaa}
\times
{}_{3}\Phi_{2}\left({-j_1\!-j_2\!+j_3,\atop
-2j_2,}{-j_1\!-j_2\!-j_3\!-1,\atop-j_1\!-j_2\!+m_3}
{-j_2\!-m_2\atop\phantom{-\lm}}
\!\!\left\vert q^2,q^{2(j_1\!+m_1\!+1)}\right.
\right).
\lb{su216}
\end{array}
\ed

We can also obtain the explicit expression of the Racah
coefficients or the $6j$-symbols for $U_q(su(2))$ using the extremal
projector.
Let $\{|j_im_i\bigr>\}$ be canonical bases of three IRs $j_i$ $(i=1,2,3)$.
The Racah coefficients for $U_q(su(2))$ (or the $q$-Racah coefficients)
are matrix elements of the transformation between two couplings
of these representations: $j_1\ot(j_2\ot j_3)$ and $(j_1\ot j_2)\ot j_3$.
i.e.
\bn
\bigr|j_1,j_2j_3(j_{23})\!:\!jm\bigr>_{\!q}=\sum_{j_{12}}
U(j_1j_2jj_3;j_{12}j_{23})_{q}\,
\bigr|j_1j_2(j_{12}),j_3)\!:\!jm\bigr>_{\!q}~,
\lb{su217}
\ed
where, for example, the vector 
$\bigr|j_1,j_2j_3(j_{23})\!:\!jm\bigr>_{\!q}$
corresponds to the first coupling scheme $j_1\ot(j_2\ot j_3)$
and it has the form:
\bn
\begin{array}{rcl}
\bigr|j_1,j_2j_3(j_{23})\!:\!jm\bigr>_{\!q}&\!\!=\!\!&
\!\!\sum\limits_{m_2m_3\atop m_1m_{23}}
\bigl(j_1m_1\,j_{23}m_{23}\bigr|jm\bigr)_{\!q}
\bigl(j_2m_2\,j_3m_3\bigr|j_{23}m_{23}\bigr)_{\!q}
\\[7pt]
&&
\phantom{possi}\times
\bigr|j_1m_1\bigr>\bigr|j_2m_2\bigr>\bigr|j_3m_3\bigr>~.
\lb{su218}
\end{array}
\ed
The $q$-Racah coefficient $U(j_1j_2jj_3;j_{12}j_{23})_{q}$ is
connected with the $q$-$6j$-symbol $\{{\ld\atop\cd}\}_q$ by
the standard relation
\bn
U(j_1j_2jj_3;j_{12}j_{23})_{q}=(-1)^{j_1\!+\!j_2\!+\!j_3\!+\!j}
\sqrt{[2j_{12}+1]![2j_{23}+1]!}
\left\{{j_1\,j_2\,j_{12}\atop j_3\,j\,j_{23}}\right\}_{\!q}.
\lb{su219}
\ed
It is not hard to obtain the formula for calculation
of the $q$-$6j$-symbols in terms of projection operators:
\bn
\!\!\!\!\!\!\!\!
\begin{array}{rcl}
&&
\!\!\!\!\!\!\!\!\!
{\ds\left\{{j_1\atop j_3}{j_2\atop j}{j_{12}\atop j_{23}}\right\}_{\!q}}=
\frac{(-1)^{j_1+j_2+j_3+j}}{\sqrt{[2j_{12}+1]![2j_{23}+1]!}}
\\[11pt]
&&
\!\!\!\!\!\!\!
\times\mbox{\ls$\frac{\langle j_1j\!-\!j_{23}|\langle j_2j_{23}\!-\!j_2|
\langle j_2j_2|P^{j_{23}}\!(23)P^{j}(123)P^{j_{12}}\!(12)
|j_1j_{12}\!-\!j_{2}\rangle|j_2j_2\rangle|j_3j\!-\!j_{12}\rangle}
{(j_1j_{12}\!-\!j_2,j_2j_2|_{12}j_{12})_{q}
(j_{12}j_{12},j_3j\!-\!j_{12}|jj)_{q}
(j_1j\!-\!j_{23},j_{23}j_{23}|jj)_{q}
(j_2j_2,j_{23}\!-\!j_2|j_{23}j_{23})_{q}}$},
\end{array}
\lb{su220}
\ed
where the notations are used: $P^j\!:=\!P_{jj}^j$,
$P^{j_{23}}(23)\!:=\!\id\ot\D_q(P^{j_{23}})$,
$P^j(123)\!:=\!(\D_q\ot\id)\D_q(P^{j})$,
$P^{j_{12}}(12)\!:=\!\D_q(P^{j_{12}})\ot\id$.

We substitute the explicit expressions for all special CGCs in
the denominator of the right side of (\ref{su220}) and then we use
the actions of the generators of $U_q(su(2))$ on the canonical basis
vectors $\big|j_im_i\bigr>$ $(i\!=\!1,2,3)$, and as result we obtain
the explicit expression for the $q$-$6j$-symbol \cite{STK2,STK3}:
\bn
\begin{array}{rcl}
&&\!\!\!\!\!\!\!\!\!\!
{\ds\left\{{j_1\atop j_3}{j_2\atop j}{j_{12}\atop j_{23}}\right\}_{\!q}}=
(-1)^{j_1+j_{12}+j_{23}+j_3}
\frac{\trid(j_1j_2j_{12})_{\!q}\trid(j_2j_3j_{23})_{\!q}}
{[j_1-j_2+j_{12}]![-j_1+j_2+j_{12}]!}
\\[9pt]
&&\phantom{possibl}
\times{\frac{\trid(j_{12}j_3j)_{\!q}\trid(j_1j_{23}j)_{\!q}
[j_{12}+j_3+j+1]![j_1+j_{23}+j+1]!}
{[j_2-j_3+j_{23}]![-j_2+j_3+j_{23}]!
[j_1-j_{23}+j]![-j_{12}+j_3+j]!}}
\\[9pt]
&&\phantom{possibl}
\times\sum\limits_{z}\frac{(-1)^z[j_1+j-j_{23}+z]![j_{3}+
j-j_{12}+z]![j_{12} +j_{23}+j_2-j-z]!}{[z]!
[j_1+j_{23}-j-z]![j_{12}+j_{3}-j-z]!
[j_2+j-j_{12}-j_{23}+z]![2j+1+z]!}~,
\end{array}
\lb{su221} 
\ed 
where we use the notation: 
\bn
\trid(j_1j_2j_3)_q=\mbox{\ls$\sqrt{\frac{[j_1+j_2-j_3]!
[j_1-j_2+j_3]![-j_1+j_2+j_3]!}{[j_1+j_2+j_3+1]!}}$}. 
\lb{su222}
\ed 
A total list of different explicit expressions and symmetry
properties for the $q$-$6j$-symbols (or the $q$-Racah
coefficients) can be found in the papers \cite{AST4,KR,STK2,STK3}.

The general formula (\ref{su221}) of the $q$-$6j$-symbol can be expressed
in terms of the following basic hypergeometric series
\bn
\begin{array}{rcl}
&&\!\!\!\!\!\!\!
{\ds\left\{{j_1\atop j_3}{j_2\atop j}{j_{12}\atop j_{23}}\right\}_{\!q}}=
(-1)^{j_1+j_{12}+j_{23}+j_3}
\frac{\trid(j_1j_2j_{12})_{\!q}\trid(j_2j_3j_{23})_{\!q}
\trid(j_{12}j_3j)_{\!q}}
{[j_1-j_2+j_{12}]![-j_1+j_2+j_{12}]![j_2-j_3+j_{23}]!}
\\[9pt]
&&
\phantom{possibl}
\times{\frac{\trid(j_1j_{23}j)_{\!q}
[j_{12}+j_3+j+1]![j_1+j_{23}+j+1]![j_2+j_{12}+j_{23}-j]!}
{[-j_2+j_3+j_{23}]![j_1+j_{23}-j]![j_{12}+j_3-j]!
[j_2+j-j_{12}-j_{23}]![2j+1]!}}
\\[9pt]
&&
\phantom{possibl}
\times{}_{4}\Phi_{3}
\left({-j_1-j_{23}+j,\,-j_{12}-j_3+j,\,j_1-j_{23}+j+1,\,j_3+j-j_{12}+1
\atop -j_2-j_{12}-j_{23}+j,\;j_2-j_{12}-j_{23}+j+1,\;2j+2}
\!\!\left\vert q^2,q^2\right.
\right).
\lb{su223}
\end{array}
\ed If we set $j_{12}\!=\!j_1\!+\!j_2$ in (\ref{su221}) we obtain
simple explicit expression for the special (so called
'stretched') $q$-$6j$-symbol: \bn
\begin{array}{rcl}
&&\!\!\!\!\!\!\!\!\!\!\!\!\!\!\!\!
{\ds\left\{{j_1\atop j_3}{j_2\atop j}
{j_1\!+\!j_2\atop j_{23}}\!\right\}_{\!q}}=
(-1)^{j_1+j_2+j_3+j}\Bigl[\frac{[2j_1]![2j_2]![j_1+j_2+j_3+j+1]!}
{[2j_1\!+2j_2\!+1]![j_1\!+j+j_{23}\!+1]!
[j_2\!+j_3\!+j_{23}\!+1]!}\Bigr]^{\frac{1}{2}}
\\[11pt]
&&\phantom{possiblea}\times
\Bigl[\frac{[j_1+j_2-j_3+j]![j_1+j_2+j_3-j]![-j_1+j+j_{23}]!
[-j_2\!+j_3\!+j_{23}]!}{[-j_1\!-j_2\!+j_3\!+j]![j_1\!+j-j_{23}]!
[j_1\!-j\!+j_{23}]![j_2\!+j_3\!-j_{23}]!
[j_2\!-j_3\!+j_{23}]!}\Bigr]^{\frac{1}{2}}.
\end{array}
\lb{su224}
\ed

In the next section we consider a more complicated example of
application of the projection operator method for calculation of
a general expression for CGCs of the quantum algebra $U_q(su(3))$.

\setcounter{equation}{0}
\section{Clebsch-Gordan coefficients for the quantum
algebra $U_q(su(3))$}
Let $\Pi:=\{\a_{1},\a_{2}\}$ be a system of
simple roots of the Lie algebra $sl(3)$
($sl(3)\!:=\!sl(3,\CC)\simeq A_2$), endowed with  the following
scalar product: $(\a_{1},\a_{1})\!=\!(\a_{2},\a_{2})\!=\!2$,
$(\a_{1},\a_{2})\!=\!(\a_{2},\a_{1})\!=\!-1$. The root system
$\Dp$ of $sl(3)$ consists of the roots $\a_1,\a_1\!+\!\a_2,\a_2$.

The quantum Hopf algebra $U_{q}(sl(3))$ is generated by the Chevalley
elements $q^{\pm h_{\a_i}}$, $e_{\pm\a_i}$ $(i=1,2)$ with the relations
(\ref{GT2}), (\ref{GT3}) where $i,j\!=\!1,2$.

For construction of the composite root vectors
$e_{\pm(\a_1+\a_2)}^{}$ we fix the following normal ordering in
$\Dp$: 
\bn 
\a_1,\;\a_1+\a_2,\;\a_2. 
\lb{su33} 
\ed
According to
this ordering we set \bn
e_{\a_1+\a_2}^{}:=[e_{\a_1}^{},e_{\a_2}^{}]_{q^{-1}}^{}~,\qquad
e_{-\a_1-\a_2}^{}:=[e_{-\a_2}^{},e_{-\a_1}^{}]_{q}^{}~. \lb{su34}
\ed Let us introduce another standard notations for the
Cartan-Weyl generators: 
\bn
\begin{array}{rcccccl}
e_{12}^{}&\!\!:=\!\!&e_{\a_1}^{},\qquad\;\;
e_{21}^{}&\!\!:=\!\!&e_{-\a_1}^{},\qquad\;\;
e_{11}^{}-e_{22}^{}&\!\!:=\!\!&h_{\a_1}^{},
\\[5pt]
e_{23}^{}&\!\!:=\!\!&e_{\a_2}^{},\qquad\;\;
e_{32}^{}&\!\!:=\!\!&e_{-\a_2}^{},\qquad\;\;
e_{22}^{}-e_{33}^{}&\!\!:=\!\!&h_{\a_2}^{},
\\[5pt]
e_{13}^{}&\!\!:=\!\!&e_{\a_1+\a_2}^{},\quad \;
e_{31}^{}&\!\!:=\!\!&e_{-\a_1-\a_2}^{},\quad\;
e_{11}^{}-e_{33}^{}&\!\!:=\!\!&h_{\a_1}^{}\!+h_{\a_2}.
\end{array}
\lb{su35}
\ed
The explicit formula for the extremal projector (\ref{EP2})
specialized to the case of $U_q(sl(3))$ has the form
\bn
p=p_{12}^{}p_{13}^{}p_{23}^{}~,
\lb{su36}
\ed
where the elements $p_{ij}$ ($1\le i<j\le3$) are given by
\bn
\begin{array}{rcl}
p_{ij}^{}&\!\!=\!\!&\sum\limits_{n=0}^{\infty}
\mbox{\ls$\frac{(-1)^n}{[n]!}$}\vph_{ij,n}^{}e_{ij}^{n}e_{ji}^{n}~,
\\[7pt]
\vph_{ij,n}^{}&\!\!=\!\!&q^{-(j-i-1)n}\Bigr\{\prod\limits_{s=1}^{n}
[e_{ii}^{}-e_{jj}^{}+j-i+s]\Bigr\}^{-1}~.
\end{array}
\lb{su37}
\ed
The extremal projector $p$ satisfies the relations:
\bn
e_{ij}^{}p=pe_{ji}^{}=0 \quad (i<j)~,
\qquad p^{2}=p~.
\lb{su38}
\ed

The quantum algebra $U_{q}(su(3))$ can be considered as the
quantum algebra $U_{q}(sl(3))$ endowed with the additional
Cartan involution $^*$:
\bn
h_{\a_i}^*=h_{\a_i}~,\qquad e_{\pm\a_i}^*=e_{\mp\a_i}^{}~,
\qquad q^*=q\;({\rm or}\;\,q^{-1})~.
\lb{su39}
\ed

Let $(\lm\mu)$ be a finite-dimensional IR of $U_q(su(3))$ with the
highest weight $(\lm\mu)$ ($\lm$ and $\mu$ are nonnegative integers).
The vector of the highest weight, denoted by the symbol
$\bigr|(\lm\mu)h\bigl>$, satisfy the relations
\bn
\begin{array}{rcl}
h_{\a_1}\bigr|(\lm\mu)h\bigl>&\!\!=\!\!&\lm\bigr|(\lm\mu)h\bigl>~,\qquad
h_{\a_2}\bigr|(\lm\mu)h\bigl>=
\mu\bigr|(\lm\mu)h\bigl>~,
\\[7pt]
e_{ij}^{}\bigr|(\lm\mu)h\bigl>&\!\!=\!\!&0\qquad (i<j)~.
\end{array}
\lb{su310} 
\ed 
Labeling of another basis vectors in the IR $(\lm\mu)$
depends on the choice of subalgebras of $U_q(u(3))$ (in other
words which reduction chain from $U_q(u(3))$ to subalgebras is
chosen). Here we use the Gelfand-Tsetlin reduction chain: 
\bn
U_q(su(3))\supset U_q(u_{Y}^{}(1))\ot U_q(su_{T}^{}(2))\supset
U_q(u_{T_{0}}^{}(1))~, \lb{su311} 
\ed 
where the subalgebra
$U_q(su_{T}^{}(2))$ is generated by the elements 
\bn
T_{+}:=e_{23}^{}~,\qquad T_{-}:=e_{32}^{}~,\qquad
T_{0}:=\mbox{\ls$\frac{1}{2}$}(e_{22}^{}-e_{33}^{})~, 
\lb{su312}
\ed 
the subalgebra $U_q(u_{T_0}^{}(1))$ is generated by
$q^{T_{0}}$, and the subalgebra $U_q(u_{Y}^{}(1))$ is generated
by $q^{Y}$ where\footnote{In the classical non-deformed case in
the elementary particle theory the subalgebra $su_{T}^{}(2)$ is
called the $T$-spin algebra and the element ${Y}$ is the
hypercharge operator.}: 
\bn
Y=-\mbox{\ls$\frac{1}{3}$}\bigl(2h_{\a_1}^{}+h_{\a_2}^{}\bigr)~.
\lb{su313} 
\ed 
In the case of the reduction chain (\ref{su311})
the basis vectors of IR $(\lm\mu)$ are denoted by 
\bn
\bigl|(\lm\mu)jtt_{z}\bigr>~. 
\lb{su314} 
\ed 
Here the quantum number set $jtt_{z}$ characterize the hypercharge 
$y$ and the $T$-spin $t$ and its projection $t_z$: 
\bn
\begin{array}{rcl}
q^{T_0}\bigl|(\lm\mu)jtt_{z}\bigr>&\!\!=\!\!&q^{t_z}
\bigl|(\lm\mu)jtt_{z}\bigr>~,
\\[7pt]
T_{\pm}\bigl|(\lm\mu)jtt_{z}\bigr>&\!\!=\!\!&\sqrt{[t\mp t_z]
[t\pm t_z\!+\!1]}\bigl|(\lm\mu)jtt_{z}\!\pm\!1\bigr>~,
\\[7pt]
q^{Y}\bigl|(\lm\mu)jtt_{z}\bigr>&\!\!=\!\!&q^{y}
\bigl|(\lm\mu)jtt_{z}\bigr>~,
\end{array}
\lb{su315}
\ed
where the parameter $j$ is connected with the eigenvalue $y$ of
the "hypercharge" operator $Y$ as follows
\bn
y=-\mbox{\ls$\frac{1}{3}$}\bigl(2\lm+\m\bigr)+2j~.
\lb{su316}
\ed
We can show (see \cite{PST1,PST2}) that the quantum numbers $jt$ 
are taken all nonnegative integers and half-integers such that
the sum $\mbox{\ls$\frac{1}{2}$}\mu+j+t$ is an integer and they 
are subjected to the constraints:
\bn
\left\{\begin{array}{rcl}
\!\!\!\!
\frac{1}{2}\mu+j-t&\!\!\ge\!\!&0~,
\\[5pt]
\!\!\!\!
\frac{1}{2}\mu-j+t&\!\!\ge\!\!&0~,
\\[5pt]
\!\!\!\!
-\frac{1}{2}\mu+j+t&\!\!\ge\!\!&0~,
\\[5pt]
\!\!\!\!
\frac{1}{2}\m+j+t&\!\!\ge\!\!&\lm+\m~.
\end{array}
\right. 
\lb{su317} 
\ed 
For every fixed $t$ the projection $t_z$ runs the values 
$t_z\!=\!-t,-t\!+\!1,\ld,t\!-\!1,t$. It is not hard to show that 
the orthonormalized vectors (\ref{su314}) can be represented in 
the following form 
\bn
\bigl\{\bigr|(\lm\mu)jtt_{z}\bigr>:=
N^{(\lm\mu)}_{\,jt}P^{\,t}_{\!t_{z};t}\,
\cR^{\,j}_{\frac{1}{2}\!\m-t}\bigr|(\lm\mu) h\bigr>\bigr\}~,
\lb{su318} 
\ed 
where $P^{\,t}_{\!t_{z};t_z'}$ is the general projection operator 
of the type (\ref{su212}) for the quantum algebra 
$U_q(su_{T}^{}(2))$, the element $\cR^{\,j}_{\frac{1}{2}\!\m-t}$ 
is a component of the irreducible tensor operator of rank $j$, 
$(-j\le j_z\le j)$: 
\bn
\cR^{\,j}_{j_z}=\mbox{\ls$\sqrt{\frac{[2j]!}{[j-j_{z}]!
[j+j_{z}]!}}$}(-1)^{j_z}q^{2j^2-j+j_z}\,e_{21}^{j-j_{z}}
e_{31}^{j+j_{z}}q^{-jh_{\alpha_1}-(j+j_z)T_0}~,
\lb{su319} 
\ed 
the normalizing factor $N^{(\lm\mu)}_{jt}$ has the form 
\bn
\begin{array}{rcl}
N^{(\lm\mu)}_{\,jt}\!\!&=\!\!&(-1)^{t-\frac{1}{2}\mu}\,
q^{(j+\frac{1}{2}\mu-t)(j+t)-2j^2+j\lm-\mu+2t}
\\[7pt]
&&{}\times
\mbox{\ls$\sqrt{
\frac{[\lm+\frac{1}{2}\m-j+t+1]![\lm+\frac{1}{2}\m-j-t]!
[\frac{1}{2}\m+j+t+1]![\frac{1}{2}\m-j+t]!}
{[\lm]![\m]![\lm+\m+1]![2j]![2t+1]!}}$}~. 
\end{array}
\lb{su320} 
\ed 
The operator 
\bn
F_{\!jtt_z;h}^{(\lm\mu)}=N^{(\lm\mu)}_{\,jt}P^{\,t}_{\!t_{z};t}
\cR^{\,j}_{\frac{1}{2}\!\m-t} 
\lb{su321} 
\ed
is called the lowering operator and the conjugate operator 
\bn
F_{\!h;jtt_z}^{(\lm\mu)}=(F_{\!jtt_z;h}^{(\lm\mu)})^* 
\lb{su322}
\ed 
is called the rising operator of $U_q(su(3))$.

If the extremal operator $p$ acts on a vector of the weight
$(\lm\mu)$ with respect to the operators $h_{\a_1}$ and $h_{\a_2}$
then we write $p^{(\lm\mu)}$. We multiply the extremal projector
$p^{(\lm\mu)}$ by the lowering and rising operators as follows
\bn
P_{jtt_z;j't't_z'}^{(\lm\mu)}=F_{\!jtt_z;h}^{(\lm\mu)}\,
p^{(\lm\mu)}F_{\!h;j't't_z'}^{(\lm\mu)}~, \lb{su323} \ed and we
assume that the resulting operator acts on a vector space of the
fixed weight
$(\!-\mbox{\ns$\frac{1}{3}$}(2\lm\!+\!\mu)\!+\!2j,t_z)$ with
respect to the operators $Y$ and $T_0$. The operator
$P_{\!jtt_z;j't't_z'}^{(\lm\mu)}$ is called the general
projection operator of $U_q(su(3))$.

For convenience we introduce the short notations: $\L:=(\lm\m)$
and $\g:=jtt_{z}$ and therefore the basis vector (\ref{su314})
will be denoted by $\bigl|\L\g\bigr>$.
Let $\{|\L_i\g_i\bigr>\}$ be Gelfand-Tsetlin bases of two IRs $\L_i$
$(i=1,2)$. Then $\{|\L_1\g_1\bigr>|\L_2\g_2\bigr>\}$ be a
uncoupled bases in the representation $\L_1\ot \L_2$ of
$U_q(su(3))\ot U_q(su(3))$. In this representation there is
another coupled basis $|\L_1\L_2\!:s\L_3\g_3\bigr>_{\!q}$ with
respect to $\D_q(U_q(su(3)))$ where the index $s$ classifies
multiple representations $\L$. We can expand the coupled basis in
terms of the uncoupled basis $\{|\L_1\g_1\bigr>|\L_2\g_2\bigr>\}$:
\bn 
\bigl|\L_1\L_2\!:s\L_3\g_3\bigr>_{\!q}=\sum_{\g_1,\g_2}^{}
\bigl(\L_1\g_1\,\L_2\g_2|s\L_3\g_3\bigr)_{\!q}\bigr|\L_1\g_1\bigr>
\bigr|\L_2\g_2\bigr>~, 
\lb{su324} 
\ed where the matrix element
$\bigl(\L_1\g_1\,\L_2\g_2|\L_3\g_3\bigr)_{\!q}$ is the
Clebsch-Gordan coefficient of $U_q(su(3))$. In just the same way
as for the non-quantized Lie algebra $su(3)$ (see
\cite{PST1,PST2}) we can show that any CGC of $U_q(su(3))$ can be
represented in terms of the linear combination of the matrix
elements of the projection operator (\ref{su323}) 
\bn
\bigl(\L_1\g_1\,\L_2\g_2|s\L_3\g_3\bigr)_{\!q}=\sum_{\g_2'}^{}
C(\g_2')\bigl<\L_1\g_1\bigr|\bigl<\L_2\g_2\bigr|
\D_q(P_{\!\g_3;h}^{\L_3})\bigr|\L_1h>\bigr|\L_2\g_2'\bigr>~.
\lb{su325} 
\ed
Classification of multiple representations $\L_3$
in the representation $\L_1\ot\L_2$ is special problem and we
shall not discuss it here. For the non-deformed algebra $su(3)$
this problem was considered in details in \cite{PST1,PST2}.
Concerning the matrix elements in the right-side of (\ref{su325})
we give here an explicit expression for the more general matrix
element: 
\bn 
\bigl<\L_1\g_1\bigr|\bigl<\L_2\g_2\bigr|
\D_q(P_{\!\g_3;\g_3'}^{\L_3})\bigr|\L_1\g_1'>\bigr|\L_2\g_2'\bigr>~.
\lb{su326} 
\ed
Using a tensor form of the projection operator (\ref{su323}) and
the Wigner-Racah calculus for the subalgebra $U_q(su(2))$ \cite{STK1} 
it is not hard to obtain the following result (see \cite{AST5}):
\bn
\!\!\!\begin{array}{lcl}
\bigl<\L_1\g_1\bigr|\bigl<\L_2\g_2\bigr|
\D_q(P_{\g_3;\g_3'}^{\L}\!)\bigr|\L_1\g_1'\bigr>\bigr|\L_2\g_2'\bigr>
\\[12pt]
=\bigl(t_1t_{1z}\,t_2t_{2z}\!\bigr|t_3t_{3z}\!\bigr)_{\!q}
\bigl(t_1t_{1z}'\,t_2t_{2z}'\!\bigr|t_3't_{3z}'\!\bigr)_{\!q}
\;\;A\!\!\!\sum\limits_{j_1'\!'j_2'\!'t_1'\!'t_2'\!'t_3'\!'}
C_{j_1'\!'j_2'\!'t_1'\!'t_2'\!'t_3'\!'}
\\[15pt]
\;\times
\left\{\!\!\!\!\!\!\!\!\begin{array}{cccc}
&j_1\!\!-\!j_1'\!'&\!\!j_2\!\!-\!j_2'\!'&
\!\!j_1\!\!+\!j_2\!\!-\!j_1'\!'\!\!-\!j_2'\!'\\
&t_1'\!'&\!\!t_2'\!'&\!\!t_3'\!'\\
&t_1&\!\!t_2&\!\!t_3
\end{array}
\!\!\!\right\}_{\!q}
\!\!\left\{\!\!\!\!\!\!\!\!\begin{array}{cccc}
&j_1'\!\!-\!j_1'\!'&\!\!j_2'\!\!-\!j_2'\!'&
\!\!j_1'\!\!+\!j_2'\!\!-\!j_1'\!'\!\!-\!j_2'\!'\\
&t_1'\!'&\!\!t_2'\!'&\!\!t_3'\!'\\
&t_1'&\!\!t_2'&\!\!t_3'
\end{array}
\!\!\!\right\}_{\!q}~.
\end{array}
\lb{su327}
\ed
Here
\bn
\!\!\!\!\begin{array}{rcl}
A&\!\!=\!\!&[\lm+1][\mu+1][\lm+\mu+2]
\\[7pt]
&&\!\!\!\!\!
\times\left[\frac{[2t_1\!+1][2t_2\!+1][2j_1\!+1]![2j_2\!+1]!
[\lm_3\!+\!\frac{1}{2}\mu_3\!-\!j_3\!+\!t_3\!+1]!
[\lm_3\!+\!\frac{1}{2}\mu_3\!-\!j_3\!-\!t_3]!}
{[\lm_1\!+\!\frac{1}{2}\mu_1\!-\!j_1\!+\!t_1\!+1]!
[\lm_1\!+\!\frac{1}{2}\mu_1\!-\!j_1\!-\!t_1]!
[\lm_2\!+\!\frac{1}{2}\mu_2\!-\!j_2\!+\!t_2\!+1]!
[\lm_2\!+\!\frac{1}{2}\mu_2\!-\!j_2\!-\!t_2]![2j_3]!}
\right]^{\frac{1}{2}}
\\[9pt]
&&\!\!\!\!\!
\times\left[\frac{[2t_1'\!+1][2t_2'\!+1][2j_1'\!+1]![2j_2'\!+1]!
[\lm_3\!+\!\frac{1}{2}\mu_3\!-\!j_3'\!+\!t_3'\!+\!1]!
[\lm_3\!+\!\frac{1}{2}\mu_3\!-\!j_3'\!-\!t_3']!}
{[\lm_1\!+\!\frac{1}{2}\mu_1\!-\!j_1'\!+\!t_1'\!+1]!
[\lm_1\!+\!\frac{1}{2}\mu_1\!-\!j_1'\!-\!t_1']!
[\lm_2\!+\!\frac{1}{2}\mu_2\!-\!j_2'\!+\!t_2'\!+1]!
[\lm_2\!+\!\frac{1}{2}\mu_2\!-\!j_2'\!-\!t_2']![2j_3']!}
\right]^{\frac{1}{2}},
\end{array}
\lb{su328}
\ed
the coefficient
$C_{j_1'\!'j_2'\!'t_1'\!'t_2'\!'t_3'\!'}$ does not contain any
`inner' summation (without summation!) and has the form 
\bn
\!\!\!\begin{array}{lcl}
C_{j_1'\!'j_2'\!'t_1'\!'t_2'\!'t_3'\!'}= 
\frac{(\!-1)^{2(j_1\!+\!j_2\!+\!j_3'\!\!-\!j_1'\!\!'\!\!-\!j_2'\!\!')}
q^{\psi}[2(j_1\!+\!j_2\!-\!j_1'\!'\!-\!j_2'\!')\!+\!1]!
[2(j_1'\!+\!j_2'\!-\!j_1'\!'\!-\!j_2'\!')\!+\!1]![2t_3'\!'\!+1]}
{[2j_1'\!']![2j_2'\!']![2j_1\!-\!2j_1'\!']![2j_2\!-\!2j_2'\!']!
[2j_1'\!-\!2j_2'\!']![2j_2'\!-\!2j_2'\!']!
[2(j_1\!+\!j_2\!-\!j_3\!-\!j_1'\!'\!-\!j_2'\!')]!}
\\[9pt]
\times\frac{[\lm_1\!+\!\frac{1}{2}\mu_1\!-\!j_1'\!'\!+\!t_1'\!'\!+\!1]!
[\lm_1\!+\!\frac{1}{2}\mu_1\!-\!j_1'\!'\!-\!t_1'\!']!
[\lm_2\!+\!\frac{1}{2}\mu_2\!-\!j_2'\!'\!+\!t_2'\!'\!+1]!
[\lm_2\!+\!\frac{1}{2}\mu_2\!-\!j_2'\!'\!-\!t_2'\!']!
[2t_1'\!'\!+1][2t_2'\!'\!+1]}{
[\lm_3\!+\!\frac{1}{2}\mu_3\!+j_1\!+\!j_2\!-\!j_3\!-\!
j_1'\!'\!-\!j_2'\!'\!+\!t_3'\!'\!+\!2]!
[\lm_3\!+\!\frac{1}{2}\mu_3\!+j_1\!+\!j_2\!-\!j_3'\!-\!
j_1'\!'\!-\!j_2'\!'\!-\!t_3'\!'\!+1]!}
\\[10pt]
\times
{\ds\left\{{j_1\!\!-\!\!j_1'\!'\atop\frac{1}{2}\mu_1}
{j_1'\!'\atop t_1}{j_1\atop t_1'\!'}\right\}_{\!q}}
{\ds\left\{{j_2\!\!-\!\!j_2'\!'\atop\frac{1}{2}\mu_2}
{j_2'\!'\atop t_2}{j_2\atop t_2'\!'}\right\}_{\!q}}
{\ds\left\{{j_3\atop t_3'\!'}
{\,{j_1\!\!+\!\!j_2\!\!-\!\!j_3\!\!-\!\!j_1'\!'\!\!-\!\!j_2'\!'\atop t_3}}
{\,{j_1\!\!+\!\!j_2\!\!-\!\!j_1'\!'\!\!-\!\!j_2'\!'\atop\frac{1}{2}\mu_3}}
\right\}_{\!q}}
\\[15pt]
\times {\ds\left\{{j_1'\!\!-\!\!j_1'\!'\atop\frac{1}{2}\mu_1}
{j_1'\!'\atop t_1'}{j_1'\atop t_1'\!'}\right\}_{\!q}}
{\ds\left\{{j_2'\!\!-\!\!j_2'\!'\atop\frac{1}{2}\mu_2}
{j_2'\!'\atop t_2'}{j_2'\atop t_2'\!'}\right\}_{\!q}}
{\ds\left\{{j_3'\atop t_3'\!'}
{\,{j_1'\!\!+\!j_2'\!\!-\!\!j_3'\!\!-\!\!j_1'\!'\!\!-\!\!j_2'\!'\atop
t_3'}}{\,{j_1'\!\!+\!\!j_2'\!\!-\!\!j_1'\!'\!\!-\!\!j_2'\!'\atop
\frac{1}{2}\mu_3}}
\right\}_{\!q}},
\end{array}
\lb{su329}
\ed
where 
\bn
\begin{array}{lcl}
\psi\!\!&=\!\!&\sum\limits_{i=1}^{2}\Bigl(2
\varphi(\lambda_i,\mu_i,j''_i,t''_i)
-\varphi(\lambda_i,\mu_i,j_i,t_i)
-\varphi(\lambda_i,\mu_i,j'_i,t'_i)
\\[7pt]
&&-\!t_i(t_i\!\!+\!1)
\!-\!t_i'(t'_i\!\!+\!1)\Big) 
\!-\!2\varphi(\lambda_3,\mu_3,j_3''\!,t_3'')
\!+\!\varphi(\lambda_3,\mu_3,j_3,t_3)
\\[7pt]
&&+\varphi(\lambda_3,\mu_3,j'_3\!,t_3')
+j_3''(4\lambda_3\!+2\mu_3\!+2)-2t_3''(t''_3\!-1)-2\mu_3
\\[7pt]
&&-(j_2\!+j_2'\!-2j''_2)(2\lambda_1\!+\mu_1\!-6j_1'')
-(j_3\!+j''_3)(j_3\!+\!j''_3\!+1)
\\[7pt]
&&-\!(j'_3\!\!+\!j''_3)(j'_3\!\!+\!j''_3\!\!+\!1)
\!+\!4(j_1\!\!-\!j_1'')(j_2\!\!-\!j_2'')\!+
\!4(j'_1\!\!-\!j_1'')(j'_2\!\!-\!j_2'')~,
\end{array}
\lb{su329'}
\ed
\bn
\begin{array}{rcl}
\varphi(\lambda,\mu,j,t)\!\!&=\!\!&\frac{1}{2}(\frac{1}{2}\mu+j-t)
(\frac{1}{2}\mu+j+t-3)+j(\lambda-2j+\!1)~,
\\[7pt]
j''_3\!\!&=\!\!&j_1\!+\!j_2\!-\!j_3\!-\!j_1''\!-\!j_2''=
j'_1\!+\!j'_2\!-\!j'_3\!-\!j_1''\!-\!j_2''.
\end{array}
\lb{su329"}
\ed

\nin
A detailed proof of the formulas (\ref{su327})-(\ref{su329"}) 
is given in \cite{AST5}.

The $q$-$9j$-symbol of $U_q(su(2))$ in (\ref{su327}) can be expressed
in terms of $q$-$6j$-symbols. This expression has the form:
\bn
\begin{array}{rcl}
\!\left\{\!\!\!\!\!\!\!\!\begin{array}{cccc}
&j_1&j_2&j_{12}\\
&j_3&j_4&j_{34}\\
&j_{13}&j_{24}&j
\end{array}
\!\!\!\right\}_{\!q}&\!\!=\!\!&\sum\limits_{z}(-1)^{2z}
q^{c(z)+c(j_{24})+c(j_{34})+c(j)}[2z+1]
\\[7pt]
&&\times
{\ds\left\{{j_{1}\atop z}{j_{2}\atop j_3}{j_{12}\atop j_{13}}\right\}_q}
{\ds\left\{{j_3\atop j}{j_4\atop j_{12}}{j_{34}\atop z}\right\}_q}
{\ds\left\{{j_{13}\atop j_4}{j_{24}\atop z}{j\atop j_2}\right\}_q},
\end{array}
\lb{su330}
\ed
where $c(j):=j(j\!+\!1)$.
In our case  $j_{12}\!=\!j_1\!+\!j_2$ and we have in the right side of
(\ref{su330}) one particular ('stretched') $q$-$6j$-symbol (\ref{su224})
which does not contain a summation. Therefore each $q$-$9j$-symbol
in the right-side of (\ref{su327}) is expressed in terms of a linear
combination of products of the general $q$-$6j$-symbols. 
Since the general $q$-$6j$-symbol
can be expressed in terms of the basic hypergeometric series
${}_{4}\Phi_{3}\left({\ld\atop\cd};q,q\right)$ therefore the
$q$-$9j$-symbol is expressed in terms of a linear combination of 
products of two basic hypergeometric series 
${}_{4}\Phi_{3}\left({\ld\atop\cd};q,q\right)$.

Thus, {\it the general matrix element (\ref{su327}) can be expressed in
terms of a linear combination of products of four basic hypergeometric
series ${}_{4}\Phi_{3}\left({\ld\atop\cd};q,q\right)$}.

\setcounter{equation}{0}
\section{Gelfand-Tsetlin basis for $U_q(u(n))$}
Let $\Pi\!:=\!\{\a_{1},\ld,\a_{n-1}\}$ be a system of simple roots
of the Lie algebra $sl(n)$ ($sl(n)\!:=\!sl(n,\CC)\simeq A_{n\!-\!1}$)
endowed with the following scalar product:
$(\a_{i},\a_{j})=(\a_{j},\a_{i})$, $(\a_{i},\a_{i})=2$,
$(\a_{i},\a_{i+1})=-1$, $(\a_i,\a_j)=0$ $((|i-j|>1)$.

The quantum Hopf algebra $U_{q}(sl(n))$ is generated by the Chevalley
elements $q^{\pm h_{\a_i}}$, $e_{\pm\a_i}$ $(i=1,2,\ld,n\!-\!1)$ with
the defining relations:
\bn
\begin{array}{rcl}
q^{h_{\a_i}}q^{-h_{\a_i}}&\!\!=\!\!&q^{-h_{\a_i}}q^{h_{\a_i}}=1~,
\qquad\;\;q^{h_{\a_i}}q^{h_{\a_j}}=q^{h_{\a_j}}q^{h_{\a_i}}~,
\\[6pt]
q^{h_{\a_i}}e_{\pm\a_j}q^{-h_{\a_i}}&\!\!=\!\!&
q^{\pm(\a_i,\a_j)}e_{\pm\a_j}~,\qquad\;
[e_{\a_i},e_{-\a_j}]=\d_{ij}\,[h_{\a_i}]~,
\\[7pt]
[e_{\pm\a_i},e_{\pm\a_j}]&\!\!=\!\!&0\qquad\qquad(|i-j|\geq 2)~,
\\[6pt]
[[e_{\pm\a_i}e_{\pm\a_j}]_{q}^{}e_{\pm\a_j}]_{q}^{}&\!\!=\!\!&0
\qquad\qquad(|i-j|=1)~,
\lb{GT2}
\end{array}
\ed
\bn
\begin{array}{rcccl}
\D_{q}(h_{\a_i})&\!\!=\!\!&h_{\a_i}\ot 1+1\ot h_{\a_i},
\qquad\qquad\qquad\;
S_{q}(h_{\a_i})&\!\!=\!\!&-h_{\a_i},
\\[7pt]
\D_{q}(e_{\pm\a_i}^{})&\!\!=\!\!&e_{\pm\a_i}^{}\ot q^{\frac{h_{\a_i}}{2}}
+q^{-\frac{h_{\a_i}}{2}}\ot e_{\pm\a_i}^{},
\qquad
S_{q}(e_{\a_i}^{})&\!\!=\!\!&-q^{\pm1}e_{\a_i}^{}.
\end{array}
\lb{GT3}
\ed
Below we shall use another basis in the Cartan subalgebra of the
algebra $sl(n)$ $(U_q(sl(n))$. Namely we set
\bn
\begin{array}{rcl}
e_{11}^{}\!\!&=&\!\!\mbox{\ls$\frac{1}{n}$}
\left((n\!-\!1)h_{\a_1}\!\!+(n\!-\!2)h_{\a_2}\!+\cd+
2h_{\a_{n-2}}\!+h_{\a_{n-1}}\!+N\right)~,
\\[6pt]
e_{22}^{}\!\!&=&\!\!\mbox{\ls$\frac{1}{n}$}
\left((n\!-\!1)h_{\a_1}\!\!+(n\!-\!2)h_{\a_2}\!+
\cd+\!2h_{\a_{n-2}}\!+h_{\a_{n-1}}\!+N\right)\!-h_{\a_1}~,
\\[2pt]
\ld&\ld&\ld\ld\ld\ld\ld\ld\ld\ld\ld\ld\ld\ld\ld\ld\ld\ld\ld\ld
\\[2pt]
e_{ii}^{}\!\!&=&\!\!\mbox{\ls$\frac{1}{n}$}
\left((n\!-\!1)h_{\a_1}\!\!+\!(n\!-\!2)h_{\a_2}\!\!+\!
\cd\!+\!2h_{\a_{n-\!2}}\!\!+\!h_{\a_{n-1}}\!\!+\!N\right)\!-\!
\!\sum\limits_{k=1}^{i-1}h_{\a_k}~,
\\[2pt]
\ld&\ld&\ld\ld\ld\ld\ld\ld\ld\ld\ld\ld\ld\ld\ld\ld\ld\ld\ld\ld
\\[2pt]
e_{nn}^{}\!\!&=&\!\!\mbox{\ls$\frac{1}{n}$}
\left(\!-h_{\a_1}\!-2h_{\a_2}\!-\cd-(n\!-\!2)h_{\a_{n-2}}\!-
(n\!-\!1)h_{\a_{n-1}}\!+N\right)~.
\end{array}
\lb{GT4}
\ed
Here $N$ is a central element of $g$ ($U_{q}(g)$), which is equal
to 0 for the case $g=sl(n)$ and $N\neq 0$ for $g=gl(n)$.
It is easy to see that
\bn
\begin{array}{rcl}
h_{\a_i}&\!\!=\!\!&e_{ii}-e_{i\!+\!1i\!+\!1}\qquad(i=1,\ld,n-1)~,
\\[6pt]
N&\!\!=\!\!&e_{11}+e_{22}+\ld+e_{nn}~.
\end{array}
\lb{GT5}
\ed
Dual elements to the ones $e_{ii}$ ($i=1,2,\ld,n$) will be denoted
by $\ep_i$ ($i=1,2,\ld, n$): $\ep_i(e_{jj})=(\ep_i,\ep_j)=\d_{ij}$.
In the terms of $\ep_i$ the positive root system $\D_{+}$ of $sl(n)$
is presented as follows
\bn
\D_{+}=\{\ep_i-\ep_j\,|\,1\le i<j\le n\}~,
\lb{GT6}
\ed
where $\ep_{i}-\ep_{i+1}$ are the simple roots:
$\a_{i}\!=\!\ep_{i}\!-\!\ep_{i+1}$ $(i\!=\!1,2,\ld,n-1)$~.
For the root vectors $e_{\ep_i\!-\!\ep_j}$ $(i\neq j)$ the another
standard notations are also used
\bn
e_{ij}:=e_{\ep_i-\ep_j}~,\qquad e_{ji}:=e_{\ep_j-\ep_i}
\qquad (1\le i<j\le n)~.
\lb{GT8}
\ed
In particular, the elements $e_{ii\!+\!1}$, $e_{i\!+\!1i}$ are the
Chevalley generators: $e_{ii\!+\!1}\!=\!e_{\a_i}$,
$e_{i\!+\!1i}\!=\!e_{-\a_i}$ ($i=1,\ld,n\!-\!1$).

For construction of the composite root vectors $e_{ij}$ $(j\neq
i\pm 1)$ we fix the following normal ordering of the positive
root system $\D_{+}$ (see \cite{T3}) \bn
(\ep_1\!-\!\ep_2),(\ep_1\!-\!\ep_3,\ep_2\!-\!\ep_3),\ld,
(\ep_1\!-\!\ep_i,\ld,\ep_{i-1}\!-\!\ep_{i}),\ld,
(\ep_1\!-\!\ep_{n},\ld,\ep_1\!-\!\ep_{n}). \lb{GT9} \ed According
to this ordering we set \bn
e_{ij}:=[e_{ik},e_{kj}]_{q^{-1}}~,\qquad
e_{ji}:=[e_{jk},e_{ki}]_{q}\qquad (1\le i<k<j\le n)~. \lb{GT10}
\ed
It should be stressed that the structure of the composite root
vectors (\ref{GT10}) is independent of the choice of the index $k$
in the r.h.s. of the definition (\ref{GT10}). In particular, one
has \bn
\begin{array}{rcrcll}
e_{ij}\!\!\!&:=&\!\!\![e_{ii\!+\!1},
e_{i\!+\!1j}\!]_{q^{\!-1}}&\!\!\!=\!\!\!&
[e_{ij\!-\!1},e_{j\!-\!1j}]_{q{\!-1}}\quad &(1\le i<j\le n)~,
\\[6pt]
e_{ji}\!\!\!&:=&\!\!\![\,e_{ji\!+1},e_{i\!+1i}]_{q}&\!\!\!=\!\!\!&
[e_{jj\!-\!1},e_{j\!-\!1i}]_{q}\quad &(1\le i<j\le n)~.
\end{array}
\lb{GT11}
\ed
The explicit formula for the extremal projector (\ref{EP2})
specialized to the case of $U_q(sl(n))$ has the form
\bn
\begin{array}{rcl}
p(U_q(sl(n))\!\!&=&\!\!p(U_q(sl(n\!-\!1))(p_{1n}^{}
p_{2n}^{}\cd p_{n-2n}^{}p_{n-1n}^{})
\\[5pt]
\!\!&=&\!\!p_{12}^{}(p_{13}^{}p_{23}^{})\cd
(p_{1i}^{}\cd p_{ii+1}^{})\cd
(p_{1n}^{}\cd p_{n-1n}^{})~,
\end{array}
\lb{GT12}
\ed
where the elements $p_{ij}$ are given by the formulas (\ref{su37})
with ($1\le i<j\le n$). The extremal projector $p:=p(U_q(sl(n))$
satisfies the eqs. (\ref{su38}) with ($1\le i<j\le n$).

The quantum algebra $U_{q}(su(n))$ can be considered as the
quantum algebra $U_{q}(sl(n))$ endowed with the additional Cartan
involution $^*$: 
\bn e_{ij}^*=e_{ji}^{}~,\qquad q^*=q\;({\rm
or}\;\,q^{-1})~. 
\lb{GT13} 
\ed 
Since the quantum algebra
$U_{q}(su(n))$ can be interpreted as the algebra $U_{q}(u(n))$
with the central element $N=0$ and the inner structure of its
representations is more easily described in terms of
$U_{q}(u(n))$, we shall consider the quantum algebra
$U_{q}(u(n))$.

Let $V^{\lm_n}$ be a finite-dimensional IR of $U_q(u(n))$ with the
highest weight $\lm_n:=(\lm_{1n},\lm_{2n},\ld,\lm_{nn})$ where
$\lm_{in}\!-\!\lm_{i+1n}$ ($i=1,\ld,n\!-\!1$) are nonnegative integers.
The vector of the highest weight, denoted by the symbol
$\bigr|\lm_n\bigl>$, satisfy the relations
\bn
\begin{array}{rcl}
q^{h_{ii}}\bigr|\lm_n\bigl>\!\!&=&\!\!
q^{\lm_{in}}\bigr|\lm_n\bigl>\qquad (1\le i\le n)~,
\\[7pt]
e_{ij}\bigr|\lm_n\bigl>\!\!&=&\!\!0\qquad\quad (i<j)~.
\end{array}
\lb{GT14} \ed
Labeling of another basis vectors in IR $V^{\lm_n}$
depends on choice of subalgebras of $U_q(u(n))$ (in other words
which reduction chain from $U_q(u(n))$ to subalgebras is chosen).
Here we use the 'so called' Gelfand-Tsetlin reduction chain: \bn
U_q(u(n))\supset U_q(u(n-1))\supset\ld\supset U_q(u(k))\supset\ld
U_q(u(1))~, \lb{GT15} \ed
where the subalgebra $U_q(u(k))$ is generated by 
$e_{ij}$ with $i,j=1,2,\ld,k$.

The following theorem can be proved.
\bnth
In the $ U_q(u(n))$-module $V^{\lm_n}$ there is the orthogonal
Gelfand-Tsetlin basis consisting of all vectors of the form
\bn
\begin{array}{rcl}
\bigl|\lm\bigr>\!\!&:=&\!\!\left|\!\!\!\!\begin{array}{llllll}
&\lm_{1n}&\lm_{2n}&\ld&\lm_{n-1n}&\lm_{nn}\\
&\lm_{1n-1}&\lm_{2n-1}&\ld&\lm_{n\!-\!1n\!-\!1}&\\
&\ld&\ld&\ld&&\\
&\lm_{12}&\lm_{22}&&&\\
&\lm_{11}&&&&\\
\end{array}\right>
\\[32pt]
\!\!&=&\!\!F_{-}(\lm_1;\lm_2)F_{-}(\lm_1;\lm_2)\cd
F_{-}(\lm_{n-1};\lm_n)\bigr|\lm_n\bigl>~,
\end{array}
\lb{GT16}
\ed
where the numbers $\lm_{ij}$ satisfy the standard inequalities
('between conditions') for the Lie algebra $u(n)$, i.e.
\bn
\lm_{ij+1}\ge\lm_{ij}\ge\lm_{i+1j+1}\qquad
{\rm for}\;\; 1\le i\le j\le n\!-\!1~.
\lb{GT17}
\ed
The lowering operators $F_{-}(\lm_k;\lm_{k+1})$, ($k=1,2\ld, n-1$),
are given by
\bneqn
F_{-}(\lm_k;\lm_{k+1})=\cN(\lm_k;\lm_{k+1})\,
p(U_q(u(k)))\prod_{i=1}^{k}(e_{k\!+\!1i})^{\lm_{ik\!+\!1}-\lm_{ik}}~,
\lb{GT18}
\\
\begin{array}{rcl}
\cN(\lm_k;\lm_{k+1})\!\!&=&\!\!
\bigg\{\prod\limits_{i=1}^{k}
\frac{[l_{ik}\!-l_{k\!+\!1k\!+\!1}\!-1]!}{[l_{ik\!+\!1}\!-l_{ik}]!
[l_{ik\!+\!1}\!-l_{k\!+\!1k\!+\!1}\!-1]!}\;\times
\\[7pt]
&&\;\;\times
\prod\limits_{1\le i<j\le k}\frac{[l_{ik\!+\!1}\!-l_{jk}]!
[l_{ik}\!-l_{jk\!+\!1}\!-1]!}{[l_{ik}\!-l_{jk}]!
[l_{ik\!+\!1}\!-l_{jk\!+\!1}\!-1]!}\bigg\}^{\frac{1}{2}}
\phantom{possi}
\end{array}
\lb{GT19}
\edeqn
where $l_{ij}\!:=\!\lm_{ij}\!-\!i$.
\edth
The explicit form of the basis vectors (\ref{GT16}) allows to
calculate the actions of the Cartan-Weyl generators on these vectors.

For the Cartan elements $q^{e_{ii}}$, ($i\!=\!1,2,\ld,n$), we
easy find that
\bn
q^{e_{ii}}\bigl|\lm\bigr>=q^{S_i-S_{i-1}}\bigl|\lm\bigr>~.
\lb{GT20} \ed Here $S_i\!=\!\sum_{j=1}^{i}\lm_{ji}$, where the
numbers $\lm_{1,i},\lm_{2i},\ld,\lm_{ii}$ are ones in the i-th
row of the pattern $\lm$ of the vector (\ref{GT16}), $S_0\!=\!0$.
The computation of the action of the generators $e_{ij}$ for
$i\neq j$ is more difficult. The procedure of this computation
goes as follows. First of all, using the explicit expression
(\ref{GT16}) for the basis vector $\bigl|\lm\bigr>$, we determine
the transformation of the basis under the action of the
generators $e_{ii\!+\!1}$ and $e_{i\!+\!1i}$ for all
$i\!=\!1,2,\ld,n\!-\!1$. Then using the inductive definition
(\ref{GT10}) for the generators $e_{ij}$ we find their action for
all $i\neq j$. The result reads as follows. 
\bnth The Cartan-Weyl generators $e_{kk\!+\!s}$ and $e_{k\!+\!sk}$
($s\!=\!1,2,\ld,n\!-\!k$) of the quantum algebra $U_q(u(n))$ act
on the Gelfand-Tsetlin basis (\ref{GT16}) according to relations:
\bneqn 
&&e_{kk\!+\!s}\bigl|\lm\bigr>=\sum_{j_i,j_2,\ld,j_s}\!\!
\bigl<\lm\!+\!\sum_{p=1}^{s}\ep_{j_p,k\!+\!p\!-\!1}\bigl|
e_{kk\!+\!s}\bigr|\lm\bigr>\,\bigl|\lm\!+\!
\sum_{p=1}^{s}\ep_{j_p,k\!+\!p\!-\!1}\bigr>, \lb{GT21}
\\
&&e_{k\!+\!sk}\bigl|\lm\bigr>=\sum_{j_i,j_2,\ld,j_s}\!\!
\bigl<\lm\!-\!\sum_{p=1}^{s}\ep_{j_p,k\!+\!p\!-\!1}\bigl|
e_{k\!+\!sk}\bigr|\lm\bigr>\,\bigl|\lm\!-\!
\sum_{p=1}^{s}\ep_{j_p,k\!+\!p\!-\!1}\bigr>.
\lb{GT22}
\edeqn
Here
\bneqn
&&\begin{array}{rcl}
\bigl<\lm\!+\!\sum\limits_{p=1}^{s}\ep_{j_p,k\!+\!p\!-\!1}\bigl|
e_{kk\!+\!s}\bigr|\lm\bigr>\!\!&=&\!\!
A_s(\lm)q^{l_{j_1k}-l_{j_sk+s-1}}
\\[7pt]
&&\times\prod\limits_{r=1}^{s}\bigl<\lm\!+\!\ep_{j_r,k\!+\!r\!-\!1}
\bigl|e_{k\!+\!r\!-\!1k\!+\!r}\bigr|\lm\bigr>,\phantom{poss}
\end{array}
\lb{GT23}
\\[10pt]
&&\begin{array}{rcl}
\bigl<\lm\!-\!\sum\limits_{p=1}^{s}\ep_{j_p,k\!+\!p\!-\!1}\bigl|
e_{k\!+\!sk}\bigr|\lm\bigr>\!\!&=&\!\!
A_{s-1}(\lm)q^{l_{j_sk+s-1}-l_{j_1k}}
\\[7pt]
&&\times\prod\limits_{r=1}^{s}\bigl<\lm\!-\!\ep_{j_r,k\!+\!r\!-\!1}
\bigl|e_{k\!+\!rk\!+\!r\!-\!1}\bigr|\lm\bigr>,\phantom{poss}
\end{array}
\lb{GT24}
\edeqn
where
\bneqn
\!\!\!\!\!\!
\begin{array}{rcl}
\bigl<\lm\!+\!\ep_{j_r,k\!+\!r\!-\!1}
\bigl|e_{k\!+\!r\!-\!1k\!+\!r}\bigr|\lm\bigr>=
\left\{-\frac{\prod\limits_{i=1}^{k+r}
[l_{ik\!+\!r}\!-l_{j_rk\!+\!r\!-\!1}]
\prod\limits_{i=1}^{k+r-2}
[l_{ik\!+\!r\!-\!2}\!-l_{j_rk\!+\!r\!-\!1}\!-1]}
{\prod\limits_{i=1\atop i\neq j_r}^{k+r-1}[l_{ik\!+\!r\!-\!1}\!-
l_{j_rk\!+\!r\!-\!1}][l_{ik\!+\!r\!-\!1}\!-
l_{j_rk\!+\!r\!-\!1}\!-\!1]}\right\}^{\frac{1}{2}}\!\!\!,
\end{array}
\lb{GT25}
\\
\!\!\!\!\!\!
\begin{array}{rcl}
\bigl<\lm\!+\!\ep_{j_r,k\!+\!r\!-\!1}
\bigl|e_{k\!+\!r\!-\!1k\!+\!r}\bigr|\lm\bigr>=
\left\{-\frac{\prod\limits_{i=1}^{k+r}
[l_{ik\!+\!r}\!-l_{j_rk\!+\!r\!-\!1}\!+\!1]
\prod\limits_{i=1}^{k+r-2}
[l_{ik\!+\!r\!-\!2}\!-l_{j_rk\!+\!r\!-\!1}]}
{\prod\limits_{i=1\atop i\neq j_r}^{k+r-1}[l_{ik\!+\!r\!-\!1}\!-
l_{j_rk\!+\!r\!-\!1}\!+1][l_{ik\!+\!r\!-\!1}\!-
l_{j_rk\!+\!r\!-\!1}]}\right\}^{\frac{1}{2}}\!\!\!,
\end{array}
\lb{GT26}
\\
\!\!\!\!\!\!
\begin{array}{rcl}
A_s(\lm)=\prod\limits_{r=1}^{s}
\frac{{\rm sign}(l_{j_{r+1}k\!+\!r}\!-l_{j_rk\!+\!r\!-\!1})}
{\{[l_{j_{r+1}k\!+\!r}\!-l_{j_rk\!+\!r\!-\!1}]
[l_{j_{r+1}k\!+\!r}\!-l_{j_rk\!+\!r\!-\!1}\!-
1]\}^{\frac{1}{2}}}\,,\phantom{poss}
\end{array}
\lb{GT27} 
\edeqn 
${\rm sign}(x)\!=\!1$  for $x\ge 0$ and ${\rm
sign}(x)\!=\!-\!1$ for $x< 0$. 
\edth 
In (\ref{GT21}), (\ref{GT22}) each summation index $j_r$ runs over integers
$1,2,\ld,k\!+\!r\!-\!1$. The symbol $\bigl|\ep_{ij}\bigr>$ 
means
the Gelfand-Tsetlin patter, which has zeros everywhere, except 1 on the
place $(ij)$.
The sum of the Gelfand-Tsetlin
patterns is given with the sums of the corresponding labels, as a
sum of matrices.

In the case $s\!=\!1$ the formulas ({\ref{GT21})-(\ref{GT27})
coincide with the results of \cite{J}, where they have been given
for the first time, however without a proof.

At the limit $q\!\to\!1$ the formulas (\ref{GT16})--(\ref{GT27})
coincides with the results of the paper \cite{AST3}.

Now we consider some generalizations of the extremal projector for
the case of the quantum algebra $U_q(sl(2,\CC))$.

\setcounter{equation}{0}
\section{`Adjoint extremal projectors' for $U_{q}(sl(2,{\bf C}))$}
Let ${\bf J}^{2}$ be the Casimir invariant for $U_{q}(sl(2,{\bf
C}))$: \bn {\bf J}^{2}=\mbox{\ls$\frac{1}{2}$}
(J_{+}J_{-}+J_{-}J_{+}+[2][J_{0}]^{2})~. \lb{AEP1} \ed One easily
verifies that: \bn {\bf J}^{2}=X+[J_{0}][J_{0}+1]~, \lb{AEP2} \ed
where we use the notation \bn X:= J_{-}J_{+}~. \lb{AEP3} \ed It
is evident that the first equation (\ref{su211}) for the extremal
projector (\ref{su29}) can be rewritten in the form \bn Xp=pX=0~.
\lb{AEP4} \ed By using the Casimir operator ${\bf J}^{2}$ we can
rewrite this equation as follows \bn {\bf J}^{2}p=p{\bf
J}^{2}=[J_{0}][J_{0}+1]p~. \lb{AEP5} \ed Thus the extremal
projector $p$ is an eigenvector for the Casimir operator ${\bf
J}^{2}$ with the eigenvalue $[J_{0}][J_{0}+1]$. The equation
(\ref{AEP4}) is similar to the algebraic equation for the
$\delta$-function: \bn x\d(x)=0~. \lb{AEP6} \ed Let us continue
this analogy. It is known that $\delta(x\!+\!\ep)$ is the
generating function for derivatives of $\delta$--function, i.e.
\bn
\d^{(n)}(x)=\mbox{\ls$\frac{d^{n}\d(x+\ep)}{d\ep^{n}}\Bigr|_{\ep=0}$}~.
\lb{AEP7} \ed Introduce an element $p(\ep)$ of the form \bn
p(\ep)=\sum_{n=-\infty}^{\infty}\mbox{\ls$\frac{(-1)^{n}}
{\overline{\Gamma}(n+1-\ep)\overline{\Gamma}(2J_{0}+2+\ep)}$}
J_{-}^{n}J_{+}^{n}~. \lb{AEP8} \ed The element $p(\ep)$ is an
analog of the generating function $\delta(x\!+\!\ep)$. Its
properties are described in the proposition. \bnpr The element
$p(\ep)$ satisfies the equation \bn
Xp(\ep)=p(\ep)X=[\ep][2J_{0}+1+\ep]p(\ep) \lb{AEP9} \ed or \bn
{\bf J}^{2}p(\ep)=p(\ep){\bf
J}^{2}=[J_{0}+\ep][J_{0}+1+\ep]p(\ep)~. \lb{AEP10} \ed \edpr {\it
Proof}. The equation (\ref{AEP9}) is easily verified by direct
calculation. The equation (\ref{AEP10}) is a consequence of
(\ref{AEP9}).

We see from eq. (\ref{AEP10}) that the element $p(\ep)$ is a eigenvector
of the Casimir operator ${\bf J}^{2}$ with the eigenvalue
$[J_{0}\!+\!\ep][J_{0}\!+\ep\!+\!1]$.

Let us introduce a scaling--derivative $\tilde{D}_{x}f(x)$ of
a function $f(x)$ depending on a variable $x$ as follows
\bn
\tl{D}_{x}f(x)=\lim_{\D x\rightarrow 0}
\mbox{\ls$\frac{f(x+\Delta x)-f(x)}{[\D x]}$}=
(q^{\frac{1}{2}}-q^{-\frac{1}{2}})(\ln q)^{-1}f_{x}^{'}(x))~,
\lb{AEP11}
\ed
where $f_{x}'$ is the usual derivative of the function $f(x)$.

Let $\tl{p}^{(n)}$, $n\!=\!0,1,2,\ld$, be scaling-derivatives of
$p(\ep)$ at $\ep\!=\!0$, i.e. \bn
\tl{p}^{(n)}=(\tl{D}_{\ep})^{n}{p(\ep)\mid}_{\ep=0}~. \lb{AEP12}
\ed \bnpr The elements $\tl{p}^{(n)}$, $n\!=\!0,1,2,\ld,$ satisfy
the algebraic equations \bn
X\tl{p}^{(n)}=\tl{p}^{(n)}X=\sum_{l=0}^{n-1}a_{l}^{n}\tl{p}^{(l)}~,
\qquad {\rm for}\;\;n=0,1,2,\ld~, \lb{AEP13} \ed where \bn
a_{l}^{n}=(q^{\frac{1}{2}}-q^{-\frac{1}{2}})^{n-l-2}
(q^{\frac{2J_{0}+1}{2}}+(-1)^{n-l} q^{-\frac{2J_{0}+1}{2}})
\mbox{\ls$\frac{n!}{l!(n-l)!}$}~. \lb{AEP14} \ed The elements
$\tl{p}^{(n)}$ can be redefined to $p^{(n)}$ such that they will
satisfy the simple equations \bn
Xp^{(n)}=p^{(n)}X=p^{(n-1)}~,\qquad {\rm for}\;\;n=0,2,1,\ld
\lb{AEP15} \ed or \bn {\bf J}^{2}p^{n}=p^{(n)} {\bf
J}^{2}=[J_{0}][J_{0}+1]p^{(n)}+p^{(n-1)}~, \qquad {\rm
for}\;\;n=0,1,2,\ld~. \lb{AEP16} \ed \edpr {\it Proof}. Applying
the scaling--differentiation operator $(\tl{D}_{\ep})^{n}$ to eq.
(\ref{AEP9}) and putting $\ep\!=\!0$ we obtain the equations
(\ref{AEP13}). Let \bn p^{n}=\sum_{l=0}^{n}d_{l}^{n}
\tilde{p}^{(l)}~. \lb{AEP17} \ed Substituting this in eq.
(\ref{AEP15}) we obtain the system of equations \bn
\sum_{k=l+1}^{n}\,d_{k}^{n}a_{l}^{k}=d_{l}^{n-1}~, \qquad{\rm
for}\;\; l=0,1,2,\ld,n-1~. \lb{AEP18} \ed This system has a
unique solution if $d_{1}^{1}\!=\![2J_{0}+1]^{-1}$ and
$d_{0}^{m}\!=\!0$ for $m\!=\!1,2,\ld,n\!-\!1$. We shall not
present the solution here, since it has a cumbersome form.

{\it Remark}. (i) The elements $p^{(n)}$, $n\!=\!1,2,\ld$, are
adjoint-vectors of the Casimir operator ${\bf J}^{2}$ with eigenvalue
$[J_{0}][J_{0}\!+\!1]$. They are joined to the eigenvector $p$ of
${\bf J}^{2}$. In this connection the element $p^{(n)}$ is called
the 'adjoint extremal projector' of the n-th order.

(ii) The elements $(-1)^{n}(n!)p^{(n)}$, $n\!=\!0,1,2,\ld$, are
analogs of the $\delta$--function and its derivatives, since they
satisfy the same algebraic equations.

(iii) In limit $q\rightarrow 1$ the elements $p^{(n)} (\tl{p}^{(n)})$,
$n=0,1,2,\ld$, turn into the corresponding elements of $sl(2,{\bf C})$
\cite{T4}.

It was found that the `adjoint extremal projectors',
$p^{(n)}(U_q(sl(2))$, $n\!=\!1,2,\ld$, are closely connected with
a special class of decomposable representations for quantum
algebra $U_q(sl(2))$ (details see \cite{DT}).

\section*{Acknowledgments}
The author is thankful to 
the NATO Advance Study Institute SF2000 for the support of his
visit. This work was supported by the Russian Foundation for
Fundamental Research, grant No. 99-01-01163, and by the program of
French-Russian scientific cooperation (CNRS grant PICS-608 and
grant RFBR-98-01-22033).

\end{document}